\documentclass[12pt]{article}

\usepackage[applemac]{inputenc}
\usepackage{amsmath,amssymb,fullpage}
\usepackage{showlabels}
\usepackage{color}

\newcommand{\dproof}{\noindent {Proof.} \quad}
\newcommand{\fproof}{\hfill $\square$ \bigskip}

\newtheorem{definition}{Definition}[section]
\newtheorem{example}{Example}[section]
\newtheorem{theorem}[definition]{Theorem}
\newtheorem{problem}[definition]{Problem}
\newtheorem{remark}[definition]{ \it Remark}
\newtheorem{coro}[definition]{Corollary}
\newtheorem{proposition}[definition]{Proposition}
\newtheorem{lemma}[definition]{Lemma}

\numberwithin{equation}{section}

\def\1B{\text{1\!\!I}}

\begin{document}
\date{13 October 2015}

\title{A Donsker delta functional approach to optimal insider control and applications to finance}

\author{
Olfa Draouil$^{1}$ and Bernt \O ksendal$^{2,3,4}$}

\footnotetext[1]{Department of Mathematics, University of Tunis El Manar, Tunis, Tunisia.
Email: {\tt olfadraouil@hotmail.fr}}

\footnotetext[2]{Department of Mathematics, University of Oslo, P.O. Box 1053 Blindern, N--0316 Oslo, Norway. Email: {\tt oksendal@math.uio.no}}

\footnotetext[3]{Norwegian School of Economics (NHH), Helleveien 30, N--5045 Bergen, Norway.
}
\footnotetext[4]{This research was carried out with support of CAS - Centre for Advanced Study, at the Norwegian Academy of Science and Letters, within the research program SEFE.}
\maketitle

\begin{abstract}
We study \emph{optimal insider control problems}, i.e. optimal control problems of stochastic systems where the controller at any time $t$, in addition to knowledge about the history of the system up to this time, also has additional information related to a \emph{future} value of the system. Since this puts the associated controlled systems outside the context of semimartingales, we apply anticipative white noise analysis, including forward integration and Hida-Malliavin calculus to study the problem. Combining this with Donsker delta functionals we transform the insider control problem into a classical (but parametrised) adapted control system, albeit with a non-classical performance functional. We establish a sufficient and a necessary maximum principle for such systems. Then we apply the results to obtain explicit solutions for some optimal insider portfolio problems in financial markets described by It\^ o-L\' evy processes.
Finally, in the Appendix we give a brief survey of the concepts and results we need from the theory of white noise, forward integrals and Hida-Malliavin calculus.
\end{abstract}



\section{Introduction}
In this paper we present a general method for solving \emph{optimal insider control problems}, i.e. optimal stochastic control problems where the controller has access to some future information about the system. This inside information in the control process puts the problem outside the context of semimartingale theory, and we therefore apply general \emph{anticipating white noise calculus}, including \emph{forward integrals} and \emph{Hida-Malliavin calculus}. Combining this with the \emph{Donsker delta functional} for the random variable $Y$ which represents the inside information, we are able to prove both a sufficient and a necessary maximum principle for the optimal control of such systems.

We then apply this machinery to the problem of optimal portfolio for an insider in a jump-diffusion financial market, and we obtain explicit expressions for the optimal insider portfolio in several cases, extending results that have been obtained earlier (by other methods) in \cite{PK}, \cite{BO}, \cite{DMOP2} and \cite{OR1}.

We now explain this in more detail:\\

The system we consider, is described by a stochastic differential equation driven by a Brownian motion $B(t)$ and an independent compensated Poisson random measure $\tilde{N}(dt,d\zeta)$, jointly defined on a filtered probability space $(\Omega, \mathbb{F}=\{ \mathcal{F}_t \}_{t \geq 0},\mathbf{P})$ satisfying the usual conditions. We assume that the inside information is of \emph{initial enlargement} type. Specifically, we assume  that the inside filtration $\mathbb{H}$ has the form

\begin{equation}\label{eq1.1}
 \mathbb{H}= \{ \mathcal{H}_t\}_{t \geq 0}, \text{ where } \mathcal{H}_t = \mathcal{F}_t \vee Y
\end{equation}
for all $t$, where $Y$ is a given $\mathcal{F}_{T_0}$-measurable random variable, for some $T_0 > T$ (both constants).
Here and in the following we choose the right-continuous version of  $ \mathbb{H}$, i.e. we put
$\mathcal{H}_{t}= \mathcal{H}_{t^+}=\bigcap_{s>t}\mathcal{H}_s.$
We assume that the value at time $t$ of our insider control process $u(t)$ is allowed to depend on both $Y$ and $\mathcal{F}_t$. In other words, $u$ is assumed to be $\mathbb{H}$-adapted. Therefore it has the form
\begin{equation}\label{eq1.2}
    u(t,\omega) = u_1(t, Y, \omega)
\end{equation}
for some function $u_1 : [0, T]\times \mathbb{R}\times\Omega\rightarrow \mathbb{R}$
such that $u_1(t, y)$ is $\mathbb{F}$-adapted for each $y \in\mathbb{R}$. For simplicity (albeit with some abuse of notation) we will in the following write $u$ in stead of $u_1$.
Consider a controlled stochastic process $X(t)=X^u(t)$ of the form
\begin{equation}\label{eq1.3}
    \left\{
  \begin{array}{l}
dX(t)=b(t,X(t),u(t),Y)dt+\sigma(t,X(t),u(t),Y)dB(t)\\
+\int_{\mathbb{R}} \gamma(t,X(t),u(t),Y,\zeta)\tilde{N}(dt,d\zeta); \quad t\geq 0\\
X(0)=x,\quad\hbox{$x\in \mathbb{R}$,}
  \end{array}
    \right.
\end{equation}
where
$u(t)=u(t,y)_{y=Y}$ is our insider control and the (anticipating) stochastic integrals are interpreted as \emph{forward integrals}, as introduced in \cite{RV} (Brownian motion case) and in \cite{DMOP1} (Poisson random measure case). A motivation for using forward integrals in the modelling of insider control is given in \cite{BO}.
We assume that the functions
\begin{align}
b(t,x,u,y)&=b(t,x,u,y,\omega):
[0,T_0]\times \mathbb{R} \times \mathbb{R} \times \mathbb{R} \times \Omega \mapsto \mathbb{R}\nonumber\\
\sigma(t,x,u,y)&=\sigma(t,x,u,y,\omega):
[0,T_0]\times \mathbb{R} \times \mathbb{R} \times \mathbb{R} \times \Omega \mapsto \mathbb{R}\nonumber\\
\gamma(t,x,u,y,\zeta)&=\gamma(t,x,u,y,\zeta,\omega):
[0,T_0]\times \mathbb{R} \times \mathbb{R} \times \mathbb{R} \times \mathbb{R} \times \Omega \mapsto \mathbb{R}\nonumber\\
\end{align}
are given bounded $C^1$ functions with respect to $x$ and $u$ and adapted processes in $(t,\omega)$ for each given $x,y,u,\zeta$, and that the forward integrals are well-defined.
Let $\mathcal{A}$ be a given family of admissible $\mathbb{H}-$adapted controls $u$.
The \emph{performance functional} $J(u)$ of a control process $u \in \mathcal{A}$ is defined by
\begin{eqnarray}\label{eq1.4}
  J(u) &=& \mathbb{E}[\int_0^T f(t, X(t),u(t),Y)dt +g(X(T),Y)],  \\ \nonumber
\end{eqnarray}
where \begin{align}
&f(t,x,u,y): [0;T] \times \mathbb{R} \times \mathbb{U} \times\mathbb{R} \mapsto \mathbb{R}\nonumber\\
&g(x,y):\mathbb{R} \times \mathbb{R} \mapsto \mathbb{R}
\end{align}
are given bounded functions, $C^1$ with respect to $x$ and $u$. The function $f$ and $g$ is called the \emph{profit rate} and \emph{terminal payoff}, respectively. For completeness of the presentation we allow these functions to depend explicitly on the future value $Y$ also, although this would not be the typical case in applications. But it could be that $f$ and $g$ are influenced by the future value $Y$ directly through the action of an insider, in addition to being influenced indirectly through the control process $u$ and the corresponding state process $x$.\\

We consider the problem to find $u^{\star} \in\mathcal{A}$ such that
\begin{equation}\label{eq1.5}
    \sup_{u\in\mathcal{A}}J(u)=J(u^{\star}).
\end{equation}
We use the Donsker delta functional of $Y$ to transform this anticipating system into a classical (albeit parametrised) adapted system with a non-classical performance functional. Then we solve this transformed system by using  modified maximum principles.\\

Here is an outline of the content of the paper:\\
\begin{itemize}
\item
In Section 2 we discuss properties of the Donsker delta functional and its conditional expectation and Hida-Malliavin derivatives.
\item
In Section 3 we present the general insider control problem and its transformation to a more classical problem.
\item
In Sections 4 and 5 we present a sufficient and a necessary maximum principle, respectively, for the transformed problem.
\item
Then in Section 6 we illustrate our results by applying them to optimal portfolio problems for an insider in a financial market.
\item
Finally, in the Appendix (Sections 7 and 8) we give a brief survey of the concepts and results we are using from white noise theory, forward integration and Hida-Malliavin calculus.
\end{itemize}

\section{The Donsker delta functional}
\begin{definition}
Let $Z :\Omega\rightarrow\mathbb{R}$ be a random variable which also belongs to the Hida space $(\mathcal{S})^{\ast}$ of stochastic distributions. Then a continuous functional
\begin{equation}\label{donsker}
    \delta_Z(.): \mathbb{R}\rightarrow (\mathcal{S})^{\ast}
\end{equation}
is called a Donsker delta functional of $Z $ if it has the property that
\begin{equation}\label{donsker property }
    \int_{\mathbb{R}}g(z)\delta_Z(z)dz = g(Z) \quad a.s.
\end{equation}
for all (measurable) $g : \mathbb{R} \rightarrow \mathbb{R}$ such that the integral converges.
\end{definition}

The Donsker delta functional is related to the \emph{regular conditional distribution}. The connection is the following:\\

As in Chapter VI in the book by Protter \cite{P}, we define the \emph{regular conditional distribution} with respect to $\mathcal{F}_t$ of a given real random variable $Y$, denoted by $Q_t(dy)=Q_t(\omega,dy)$, by the following properties:
\begin{itemize}
\item
For any Borel set $\Lambda \subseteq \mathbb{R}, Q_t(\cdot, \Lambda)$ is a version of $\mathbb{E}[\mathbf{1}_{Y \in \Lambda} | \mathcal{F}_t]$
\item
For each fixed $\omega, Q_t(\omega,dy)$ is a probability measure on the Borel subsets of $\mathbb{R}$
\end{itemize}

It is well-known that such a regular conditional distribution always exists. See e. g. \cite{B}, page 79.

From the required properties of $Q_t(\omega,dy)$ we get the following formula
\begin{equation}
\int_{\mathbb{R}} f(y) Q_t(\omega,dy) = \mathbb{E}[ f(Y) | \mathcal{F}_t]
\end{equation}
Comparing with the definition of the Donsker delta functional, we obtain the following representation of the regular conditional distribution:

\begin{proposition}
Suppose $Q_t(\omega,dy)$ is absolutely continuous with respect to Lebesgue measure on $\mathbb{R}$. Then the Donsker delta functional of $Y$,  $\delta_Y(y)$, exists and we have
\begin{equation}
\frac{Q_t(\omega,dy)}{dy} = \mathbb{E}[ \delta_Y(y) | \mathcal{F}_t]
\end{equation}
\end{proposition}

A general expression, in terms of Wick calculus, for the Donsker delta functional of an It\^ o diffusion with non-degenerate diffusion coefficient can be found in the amazing paper \cite{LP}. See also \cite{MP}. In the following we present more explicit formulas the Donsker delta functional and its conditional expectation and Hida-Malliavin derivatives, for It\^{o}-L\'{e}vy processes:

\subsection{The Donsker delta functional for a class of It\^{o} - L\'{e}vy processes}
Consider the special case when $Y$ is a first order chaos random variable of the form
\begin{equation}\label{eq2.5}
    Y = Y (T_0); \text{ where } Y (t) =\int_0^t\beta(s)dB(s)+\int_0^t\int_{\mathbb{R}}\psi(s,\zeta)\tilde{N}(ds,d\zeta), \mbox{ for } t\in [0,T_0]
\end{equation}
for some deterministic functions $\beta \neq 0, \psi$ satisfying
\begin{equation}\label{}
    \int_0^{T_0} \{ \beta^2(t)+\int_{\mathbb{R}}\psi^2(t,\zeta)\nu(d\zeta)\} dt<\infty \text{ a.s. }
\end{equation}
We also assume that the growth condition \eqref{eq8.4} holds throughout this paper.\\

In this case it is well known (see e.g. \cite{MOP}, \cite{DO}, Theorem 3.5, and \cite{DOP},\cite{DiO}) that the Donsker delta functional exists in $(\mathcal{S})^{\ast}$ and is given
by
\begin{eqnarray}\label{eq2.7}
   \delta_Y(y)&=&\frac{1}{2\pi}\int_{\mathbb{R}}\exp^{\diamond}\big[ \int_0^{T_0}\int_{\mathbb{R}}(e^{ix\psi(s,\zeta)}-1)\tilde{N}(ds,d\zeta)+ \int_0^{T_0}ix\beta(s)dB(s)  \nonumber\\
   &+&  \int_0^{T_0}\{\int_{\mathbb{R}}(e^{ix\psi(s,\zeta)}-1-ix\psi(s,\zeta))\nu(d\zeta)-\frac{1}{2}x^2\beta^2(s)\}ds-ixy\big]dx.
\end{eqnarray}

We will need an expression for the conditional expectation
$$ \mathbb{E}[\delta_Y(y)|\mathcal{F}_t].$$
To this end, we proceed as follows:

Using the Wick rule when taking conditional expectation, using the martingale properties of
the processes $\int_0^t\int_{\mathbb{R}}(e^{ix\psi(s,\zeta)}-1)\tilde{N}(ds,d\zeta)$ and $\int_0^t \beta(s)dB(s)$, we get:

\begin{eqnarray}\label{eq2.9}
  \mathbb{E}[\delta_Y(y)|\mathcal{F}_t] &=&\frac{1}{2\pi}\int_{\mathbb{R}}\mathbb{E}\big[\exp^{\diamond}\big[ \int_0^{T_0}\int_{\mathbb{R}}(e^{ix\psi(s,\zeta)}-1)\tilde{N}(ds,d\zeta)+ \int_0^{T_0}ix\beta(s)dB(s)  \nonumber\\
  &+& \int_0^{T_0}\{\int_{\mathbb{R}}(e^{ix\psi(s,\zeta)}-1-ix\psi(s,\zeta))\nu(d\zeta)-\frac{1}{2}x^2\beta^2(s)\}ds-ixy\big]|\mathcal{F}_t\big]dx \nonumber\\
   &=&\frac{1}{2\pi} \int_{\mathbb{R}} \exp^{\diamond}\big[ \mathbb{E}\big[ \int_0^{T_0}\int_{\mathbb{R}}(e^{ix\psi(s,\zeta)}-1)\tilde{N}(ds,d\zeta)+ \int_0^{T_0}ix\beta(s)dB(s)  \nonumber\\
   &+& \int_0^{T_0}\{\int_{\mathbb{R}}(e^{ix\psi(s,\zeta)}-1-ix\psi(s,\zeta))\nu(d\zeta)-\frac{1}{2}x^2\beta^2(s)\}ds-ixy |\mathcal{F}_t\big]\big]dx  \nonumber\\
   &=&\frac{1}{2\pi}\int_{\mathbb{R}}\exp^{\diamond}\big[ \int_0^t\int_{\mathbb{R}}(e^{ix\psi(s,\zeta)}-1)\tilde{N}(ds,d\zeta)+ \int_0^tix\beta(s)dB(s)    \nonumber\\
   &+& \int_0^{T_0}\{\int_{\mathbb{R}}(e^{ix\psi(s,\zeta)}-1-ix\psi(s,\zeta))\nu(d\zeta)-\frac{1}{2}x^2\beta^2(s)\}ds-ixy\big]dx \nonumber\\
   &=& \frac{1}{2\pi}\int_{\mathbb{R}}\{\exp^{\diamond}\big[ \int_0^t\int_{\mathbb{R}}(e^{ix\psi(s,\zeta)}-1)\tilde{N}(ds,d\zeta)\big]\}\diamond \{\exp^{\diamond}\big[ \int_0^t ix\beta(s)dB(s)\big]\}  \nonumber\\
   &\diamond& \{\exp^{\diamond}\big[\int_0^{T_0}\{\int_{\mathbb{R}}(e^{ix\psi(s,\zeta)}-1-ix\psi(s,\zeta))\nu(d\zeta)-\frac{1}{2}x^2\beta^2(s)\}ds-ixy\big]\}dx \nonumber\\
   &=&  \frac{1}{2\pi}\int_{\mathbb{R}}\exp\big[\int_0^t\int_{\mathbb{R}}ix\psi(s,\zeta)\tilde{N}(ds,d\zeta) +\int_0^t ix\beta(s)dB(s)\\
   &+&\int_t^{T_0}\int_{\mathbb{R}}(e^{ix\psi(s,\zeta)}-1-ix\psi(s,\zeta))\nu(d\zeta)ds-\int_t^{T_0}\frac{1}{2}x^2\beta^2(s)ds-ixy\big]dx\nonumber
\end{eqnarray}

Here we have used that (see e.g. \cite{DO},Lemma 3.1)
\begin{equation}
\exp^{\diamond}\big[ \int_0^{T_0}ix\beta(s)dB(s)\big]=\exp\big[\int_0^{T_0}ix\beta(s) dB(s) +\frac{1}{2}\int_0^{T_0}x^2 \beta^2(s) ds \big]
 \end{equation}
 and
 \begin{align}
&\exp^{\diamond}\big[ \int_0^{T_0}\int_{\mathbb{R}}(e^{ix\psi(s,\zeta)}-1)\tilde{N}(ds,d\zeta)]\nonumber\\
&= \exp\big[\int_0^{T_0}\int_{\mathbb{R}}ix\psi(s,\zeta)\tilde{N}(ds,d\zeta)- \int_0^{T_0}\int_{\mathbb{R}}(e^{ix\psi(s,\zeta)}-1-ix\psi(s,\zeta))\nu(d\zeta)  \big]
\end{align}

We proceed to find
$$ \mathbb{E}[D_{t,z}\delta_Y(y)| \mathcal{F}_t],$$
where $D_{t,\zeta}$ denotes the Hida-Malliavin derivative at $(t,\zeta) \in [0,T] \times \mathbb{R}$ with respect to the Poisson random measure $N$:\\
First, note that
\begin{eqnarray}\label{eq5}
 D_{t,z}\delta_Y(y)&=&\frac{1}{2\pi}\int_{\mathbb{R}} D_{t,z}\exp^{\diamond}\big[ \int_0^{T_0}\int_{\mathbb{R}}(e^{ix\psi(s,\zeta)}-1)\tilde{N}(ds,d\zeta)+ \int_0^{T_0}ix\beta(s)dB(s)  \nonumber\\
   &+&\int_0^{T_0}\{\int_{\mathbb{R}}(e^{ix\psi(s,\zeta)}-1-ix\psi(s,\zeta))\nu(d\zeta)-\frac{1}{2}x^2\beta^2(s)\}ds-ixy\big]dx  \nonumber\\
   &=&\frac{1}{2\pi}\int_{\mathbb{R}}\exp^{\diamond}\big[ \int_0^{T_0}\int_{\mathbb{R}}(e^{ix\psi(s,\zeta)}-1)\tilde{N}(ds,d\zeta)+ \int_0^{T_0}ix\beta(s)dB(s)  \nonumber\\
   &+&\int_0^{T_0}\{\int_{\mathbb{R}}(e^{ix\psi(s,\zeta)}-1-ix\psi(s,\zeta))\nu(d\zeta)-\frac{1}{2}x^2\beta^2(s)\}ds-ixy\big]  \nonumber\\
   &\diamond&D_{t,z} \big[ \int_0^{T_0}\int_{\mathbb{R}}(e^{ix\psi(s,\zeta)}-1)\tilde{N}(ds,d\zeta)+ \int_0^{T_0}ix\beta(s)dB(s) \nonumber\\
   &+&\int_0^{T_0}\{\int_{\mathbb{R}}(e^{ix\psi(s,\zeta)}-1-ix\psi(s,\zeta))\nu(d\zeta)-\frac{1}{2}x^2\beta^2(s)\}ds-ixy\big]dx \nonumber\\
   &=&\frac{1}{2\pi}\int_{\mathbb{R}}\exp^{\diamond}\big[ \int_0^{T_0}\int_{\mathbb{R}}(e^{ix\psi(s,\zeta)}-1)\tilde{N}(ds,d\zeta)+ \int_0^{T_0}ix\beta(s)dB(s)  \nonumber\\
   &+& \int_0^{T_0}\{\int_{\mathbb{R}}(e^{ix\psi(s,\zeta)}-1-ix\psi(s,\zeta))\nu(d\zeta)-\frac{1}{2}x^2\beta^2(s)\}ds-ixy\big]  \nonumber\\
   &\times&(e^{ix\psi(t,z)}-1)dx
\end{eqnarray}
Here we have used that
 $$D_{t,z} \int_0^T \beta(s) dB(s) = 0,$$
 which follows from our assumption that $B$ and $\tilde{N}$ are independent, so for $D_{t,z}$ the random variable $B(s)$ is like a constant.\\
Using equation (\ref{eq5}) and the Wick chain rule we get
\begin{eqnarray}\label{eq2.13}
 \mathbb{E}[D_{t,z}\delta_Y(y)|\mathcal{F}_t]&=& \frac{1}{2\pi}\int_{\mathbb{R}}\mathbb{E}\big[\exp^{\diamond}\big[ \int_0^{T_0}\int_{\mathbb{R}}(e^{ix\psi(s,\zeta)}-1)\tilde{N}(ds,d\zeta)+ \int_0^{T_0}ix\beta(s)dB(s)  \nonumber\\
   &+& \int_0^{T_0}\{\int_{\mathbb{R}}(e^{ix\psi(s,\zeta)}-1-ix\psi(s,\zeta))\nu(d\zeta)-\frac{1}{2}x^2\beta^2(s)\}ds-ixy\big] \nonumber\\
   &\times&(e^{ix\psi(t,z)}-1)dx|\mathcal{F}_t\big]dx \nonumber\\
   &=&\frac{1}{2\pi}\int_{\mathbb{R}}\mathbb{E}\big[\exp^{\diamond}\big[ \int_0^{T_0}\int_{\mathbb{R}}(e^{ix\psi(s,\zeta)}-1)\tilde{N}(ds,d\zeta)+ \int_0^{T_0}ix\beta(s)dB(s)  \nonumber\\
   &+& \int_0^{T_0}\{\int_{\mathbb{R}}(e^{ix\psi(s,\zeta)}-1-ix\psi(s,\zeta))\nu(d\zeta)-\frac{1}{2}x^2\beta^2(s)\}ds-ixy\big]|\mathcal{F}_t\big] \nonumber\\
   &\times& (e^{ix\psi(t,z)}-1)dx\nonumber\\
   &=& \frac{1}{2\pi}\int_{\mathbb{R}}\exp\big[\int_0^t\int_{\mathbb{R}}ix\psi(s,\zeta)\tilde{N}(ds,d\zeta) +\int_0^t ix\beta(s)dB(s)\nonumber\\
   &+&\int_t^{T_0}\int_{\mathbb{R}}(e^{ix\psi(s,\zeta)}-1-ix\psi(s,\zeta))\nu(d\zeta)ds-\int_t^{T_0}\frac{1}{2}x^2\beta^2(s)ds-ixy\big]\nonumber\\
   &\times&(e^{ix\psi(t,z)}-1)dx.
 \end{eqnarray}
 Next we want to find
$$ \mathbb{E}[D_{t}\delta_Y(y)| \mathcal{F}_t],$$
where $D_t$ denotes the Hida-Malliavin drivative at $t$ with respect to Brownian motion $B$:\\
Note that:
 \begin{eqnarray}\label{eq2.14}
 D_{t}\delta_Y(y)&=&\frac{1}{2\pi}\int_{\mathbb{R}} D_{t}\exp^{\diamond}\big[ \int_0^{T_0}\int_{\mathbb{R}}(e^{ix\psi(s,\zeta)}-1)\tilde{N}(ds,d\zeta)+ \int_0^{T_0}ix\beta(s)dB(s)  \nonumber\\
   &+&\int_0^{T_0}\{\int_{\mathbb{R}}(e^{ix\psi(s,\zeta)}-1-ix\psi(s,\zeta))\nu(d\zeta)-\frac{1}{2}x^2\beta^2(s)\}ds-ixy\big]dx  \nonumber\\
   &=&\frac{1}{2\pi}\int_{\mathbb{R}}\exp^{\diamond}\big[ \int_0^{T_0}\int_{\mathbb{R}}(e^{ix\psi(s,\zeta)}-1)\tilde{N}(ds,d\zeta)+ \int_0^{T_0}ix\beta(s)dB(s)  \nonumber\\
   &+&\int_0^{T_0}\{\int_{\mathbb{R}}(e^{ix\psi(s,\zeta)}-1-ix\psi(s,\zeta))\nu(d\zeta)-\frac{1}{2}x^2\beta^2(s)\}ds-ixy\big]  \nonumber\\
   &\diamond&D_{t} \big[ \int_0^{T_0}\int_{\mathbb{R}}(e^{ix\psi(s,\zeta)}-1)\tilde{N}(ds,d\zeta)+ \int_0^{T_0}ix\beta(s)dB(s) \nonumber\\
   &+&\int_0^{T_0}\{\int_{\mathbb{R}}(e^{ix\psi(s,\zeta)}-1-ix\psi(s,\zeta))\nu(d\zeta)-\frac{1}{2}x^2\beta^2(s)\}ds-ixy\big]dx \nonumber\\
   &=&\frac{1}{2\pi}\int_{\mathbb{R}}\exp^{\diamond}\big[ \int_0^{T_0}\int_{\mathbb{R}}(e^{ix\psi(s,\zeta)}-1)\tilde{N}(ds,d\zeta)+ \int_0^{T_0}ix\beta(s)dB(s)  \nonumber\\
   &+& \int_0^{T_0}\{\int_{\mathbb{R}}(e^{ix\psi(s,\zeta)}-1-ix\psi(s,\zeta))\nu(d\zeta)-\frac{1}{2}x^2\beta^2(s)\}ds-ixy\big]  \nonumber\\
   &\times& ix\beta(t)dx
\end{eqnarray}
Here we have used that
$$D_{t} \int_0^{T_0}\int_{\mathbb{R}}(e^{ix\psi(s,\zeta)}-1)\tilde{N}(ds,d\zeta) =0,$$
 which follows from the assumption that $B$ and $\tilde{N}$ are independent, so for $D_{t}$ the random variable $\tilde{N}(s,\zeta)$ is like a constant.\\
Using equation (\ref{eq2.14}) and the Wick chain rule we get
\begin{eqnarray}\label{2.15}
 \mathbb{E}[D_{t}\delta_Y(y)|\mathcal{F}_t]&=& \frac{1}{2\pi}\int_{\mathbb{R}}\mathbb{E}\big[\exp^{\diamond}\big[ \int_0^{T_0}\int_{\mathbb{R}}(e^{ix\psi(s,\zeta)}-1)\tilde{N}(ds,d\zeta)+ \int_0^{T_0}ix\beta(s)dB(s)  \nonumber\\
   &+& \int_0^{T_0}\{\int_{\mathbb{R}}(e^{ix\psi(s,\zeta)}-1-ix\psi(s,\zeta))\nu(d\zeta)-\frac{1}{2}x^2\beta^2(s)\}ds-ixy\big] \nonumber\\
   &\times& ix\beta(t)dx|\mathcal{F}_t\big]dx \nonumber\\
   &=&\frac{1}{2\pi}\int_{\mathbb{R}}\mathbb{E}\big[\exp^{\diamond}\big[ \int_0^{T_0}\int_{\mathbb{R}}(e^{ix\psi(s,\zeta)}-1)\tilde{N}(ds,d\zeta)+ \int_0^{T_0}ix\beta(s)dB(s)  \nonumber\\
   &+& \int_0^{T_0}\{\int_{\mathbb{R}}(e^{ix\psi(s,\zeta)}-1-ix\psi(s,\zeta))\nu(d\zeta)-\frac{1}{2}x^2\beta^2(s)\}ds-ixy\big]|\mathcal{F}_t\big] \nonumber\\
   &\times& ix\beta(t)dx\nonumber\\
   &=& \frac{1}{2\pi}\int_{\mathbb{R}}\exp\big[\int_0^t\int_{\mathbb{R}}ix\psi(s,\zeta)\tilde{N}(ds,d\zeta) +\int_0^t ix\beta(s)dB(s)\nonumber\\
   &+&\int_t^{T_0}\int_{\mathbb{R}}(e^{ix\psi(s,\zeta)}-1-ix\psi(s,\zeta))\nu(d\zeta)ds-\int_t^{T_0}\frac{1}{2}x^2\beta^2(s)ds-ixy\big]\nonumber\\
   &\times&ix\beta(t)dx.
 \end{eqnarray}

 \subsection{The Donsker delta functional for a Gaussian process}

 Consider the special case when $Y$ is a Gaussian random variable of the form
\begin{equation}\label{eq5.47}
    Y = Y (T_0); \text{ where } Y (t) =\int_0^t\beta(s)dB(s), \mbox{ for } t\in [0,T_0]
\end{equation}
for some deterministic function $\beta\in \mathbf{L}^2[0,T_0]$ with
\begin{equation}\label{}
    \|\beta\|^2_{[t,T]} :=\int_t^T\beta(s)^2ds>0 \mbox{ for all } t\in[0,T].
\end{equation}
In this case it is well known that the Donsker delta functional is given
by
\begin{equation}\label{}
    \delta_{Y}(y)=(2\pi v)^{-\frac{1}{2}}\exp^{\diamond}[-\frac{(Y-y)^{\diamond2}}{2v}]
\end{equation}
where we have put $v :=\|\beta\|^2_{[0,T_0]}$. See e.g. \cite{AaOU}, Proposition $3.2$.
Using the Wick rule when taking conditional expectation, using the martingale property of
the process $Y (t)$ and applying Lemma $3.7$ in \cite{AaOU} we get
\begin{eqnarray}\label{eq5.50}
   \mathbb{E}[\delta_Y(y)|\mathcal{F}_t]&=&(2\pi v)^{-\frac{1}{2}}\exp^{\diamond}[-\mathbb{E}[\frac{(Y(T_0)-y)^{\diamond 2}}{2v}|\mathcal{F}_t]] \nonumber \\
   &=& (2\pi \|\beta\|^2_{[0,T_0]})^{-\frac{1}{2}}\exp^{\diamond}[- \frac{(Y(t)-y)^{\diamond 2}}{2\|\beta\|^2_{[0,T_0]}}]\nonumber \\
   &=& (2\pi \|\beta\|^2_{[t,T_0]})^{-\frac{1}{2}} \exp[- \frac{(Y(t)-y)^2}{2\|\beta\|^2_{[t,T_0]}}].
\end{eqnarray}
Similarly, by the Wick chain rule and Lemma $3.8$ in \cite{AaOU} we get, for $t \in [0,T],$
\begin{eqnarray}\label{eq5.51}
  \mathbb{E}[D_t\delta_Y(y)|\mathcal{F}_t] &=&-\mathbb{E}[(2\pi v)^{-\frac{1}{2}}\exp^{\diamond}[- \frac{(Y(T_0)-y)^{\diamond 2}}{2v}]\diamond\frac{Y(T_0)-y}{v}\beta(t)|\mathcal{F}_t] \nonumber\\
   &=&-(2\pi v)^{-\frac{1}{2}} \exp^{\diamond}[- \frac{(Y(t)-y)^{\diamond 2}}{2v}]\diamond\frac{Y(t)-y}{v}\beta(t)\nonumber \\
   &=& -(2\pi \|\beta\|^2_{[t,T_0]})^{-\frac{1}{2}}\exp[- \frac{(Y(t)-y)^2}{2\|\beta\|^2_{[t,T_0]}}]\frac{Y(t)-y}{\|\beta\|^2_{[t,T_0]}}\beta(t).
\end{eqnarray}

\subsection{The Donsker delta functional for a Brownian-Poisson process}
Next, assume that $Y=Y(T_0)$, with
\begin{equation}\label{Brow-Poiss}
 Y(t)=\beta B(t)+\tilde{N}(t); \quad 0 \leq t \leq T_0
 \end{equation}
where $\beta \neq 0$ is a constant. Here $\tilde{N}(t)=N(t)-\lambda t$, where $N(t)$ is a Poisson process with intensity $\lambda > 0$.
In this case the
L\'{e}vy measure is $\nu(d\zeta) = \lambda\delta_1(d\zeta)$ since the jumps are of size $1$.
Comparing with \eqref{eq2.7} and by taking
 $\psi=1$, we obtain
\begin{equation}\label{}
    \delta_Y(y)=\frac{1}{2\pi}\int_{\mathbb{R}}\exp^{\diamond}\big[(e^{ix}-1)\tilde{N}(T_0)+i x \beta B(T_0)+\lambda T_0 (e^{ix}-1-ix)-\frac{1}{2}x^2 \beta^2 T_0-ixy\big]dx
\end{equation}

By using the general expressions \eqref{eq2.9} and \eqref{eq2.13} in Section 2.1, we get:
\begin{equation}\label{eq2.22}
    \mathbb{E}[\delta_Y(y)|\mathcal{F}_t]=\int_{\mathbb{R}}F(t,x)dx,
\end{equation}
where
\begin{equation}
F(t,x)=\frac{1}{2\pi}\exp\big[ ix\tilde{N}(t)+i x \beta B(t)+\lambda(T_0-t)(e^{ix}-1-ix)-\frac{1}{2}x^2 \beta^2 (T_0-t)-ixy\big].
\end{equation}
This gives
\begin{equation}\label{eq2.22a}
  \mathbb{E}[D_{t}\delta_Y(y)|\mathcal{F}_t]=\int_{\mathbb{R}}F(t,x)ix\beta dx
  \end{equation}
  and
\begin{equation}\label{eq2.23}
    \mathbb{E}[D_{t,1}\delta_Y(y)|\mathcal{F}_t]=\int_{\mathbb{R}}F(t,x)(e^{ix}-1)dx.
\end{equation}

\section{The general insider optimal control problem}

We now present a general method, based on the Donsker delta functional, for solving optimal insider control problems when the inside information is of \emph{initial enlargement} type. Specifically, let us from now on assume  that the inside filtration $\mathbb{H}$ has the form

\begin{equation}\label{H_t}
 \mathbb{H}= \{ \mathcal{H}_t\}_{t \geq 0}, \text{ where } \mathcal{H}_t = \mathcal{F}_t \vee Y
\end{equation}
for all $t$, where $Y \in L^2(P)$ is a given $\mathcal{F}_{T_0}$-measurable random variable, for some $T_0 > T$. We also assume that $Y$ has a Donsker delta functional $\delta_Y(y) \in (\mathcal{S})^{\ast}.$ We consider the situation when the value at time $t$ of our insider control process $u(t)$ is allowed to depend on both $Y$ and $\mathcal{F}_t$. In other words, $u$ is assumed to be $\mathbb{H}$-adapted. Therefore it has the form
\begin{equation}\label{u(t)}
    u(t,\omega) = u_1(t, Y, \omega)
\end{equation}
for some function $u_1 : [0, T]\times \mathbb{R}\times\Omega\rightarrow \mathbb{R}$
such that $u_1(t, y)$ is $\mathbb{F}$-adapted for each $y \in\mathbb{R}$. For simplicity (albeit with some abuse of notation) we will in the following write $u$ in stead of $u_1$.
Consider a controlled stochastic process $X(t)=X^u(t)$ of the form
\begin{equation}\label{richesse}
    \left\{
  \begin{array}{l}
dX(t)=b(t,X(t),u(t),Y)dt+\sigma(t,X(t),u(t),Y)dB(t)\\
+\int_{\mathbb{R}} \gamma(t,X(t),u(t),Y,\zeta)\tilde{N}(dt,d\zeta); \quad t\geq 0\\
X(0)=x,\quad\hbox{$x\in \mathbb{R}$,}
  \end{array}
    \right.
\end{equation}
with coefficients as in \eqref{eq1.3}, and where
$u(t)=u(t,y)_{y=Y}$ is our insider control.  As pointed out in the Introduction we interpret the stochastic integrals as forward integrals.\\
Then $X(t)$ is $\mathbb{H}$-adapted, and hence using the definition of the Donsker delta functional $\delta_Y(y)$ of $Y$ we get
\begin{equation}\label{eq6}
X(t)=x(t,Y)=x(t,y)_{y=Y}=\int_{\mathbb{R}}x(t,y)\delta_Y(y)dy
\end{equation}
for some $y$-parametrized process $x(t,y)$ which is $\mathbb{F}$-adapted for each $y$.
Then, again by the definition of the Donsker delta functional and the properties of forward integration (see Lemma 7.20 and Lemma 8.12), we can write
\begin{align}\label{eq7}
X(t)&= x +\int_0^t b(s,X(s),u(s),Y)ds + \int_0^t \sigma(s,X(s),u(s),Y)dB(s)\nonumber\\
&+\int_0^t \int_{\mathbb{R}} \gamma(s,X(s),u(s),Y,\zeta)\tilde{N}(ds,d\zeta)\nonumber\\
&= x+\int_0^t b(s,x(s,Y),u(s,Y),Y)ds + \int_0^t \sigma(s,x(s,Y),u(s,Y),Y)dB(s)\nonumber\\
&+\int_0^t \int_{\mathbb{R}} \gamma(s,x(s,Y),u(s,Y),Y,\zeta)\tilde{N}(ds,d\zeta)\nonumber\\
&= x+\int_0^t b(s,x(s,y),u(s,y),y)_{y=Y}ds + \int_0^t \sigma(s,x(s,y),u(s,y),y)_{y=Y}dB(s)\nonumber\\
&+\int_0^t \int_{\mathbb{R}} \gamma(s,x(s,y),u(s,y),y,\zeta)_{y=Y}\tilde{N}(ds,d\zeta)\nonumber\\
&= x+\int_0^t \int_{\mathbb{R}}b(s,x(s,y),u(s,y),y)\delta_Y(y)dyds + \int_0^t \int_{\mathbb{R}}\sigma(s,x(s,y),u(s,y),y)\delta_Y(y)dydB(s)\nonumber\\
&+\int_0^t \int_{\mathbb{R}}\int_{\mathbb{R}} \gamma(s,x(s,y),u(s,y),y,\zeta)\delta_Y(y)\tilde{N}(ds,d\zeta)\nonumber\\
&= x+\int_{\mathbb{R}} [\int_0^t b(s,x(s,y),u(s,y),y)ds + \int_0^t \sigma(s,x(s,y),u(s,y),y)dB(s)\nonumber\\
&+\int_0^t \int_{\mathbb{R}} \gamma(s,x(s,y),u(s,y),y,\zeta)\tilde{N}(ds,d\zeta)]\delta_Y(y)dy
\end{align}
Comparing \eqref{eq6} and \eqref{eq7} we see that  \eqref{eq6} holds if we choose $x(t,y)$ for each $y $ as the solution of the classical SDE
\begin{equation}\label{eq8}
    \left\{
\begin{array}{l}
    dx(t,y) = b(t,x(t,y),u(t,y),y)dt + \sigma(t,x(t,y),u(t,y),y)dB(t)\\
    + \int_{\mathbb{R}} \gamma(t,x(t,y),u(t,y),y,\zeta)\tilde{N}(dt,d\zeta); \quad t\geq 0\\
    x(0,y)  = x,\quad\hbox{$x\in \mathbb{R}$,}
\end{array}
    \right.
\end{equation}
Let $\mathcal{A}$ be a given family of admissible $\mathbb{H}-$adapted controls $u$.
The \emph{performance functional} $J(u)$ of a control process $u \in \mathcal{A}$ is defined by
\begin{eqnarray}\label{performance}
  J(u) &=& \mathbb{E}[\int_0^T f(t, X(t),u(t))dt +g(X(T))]  \\ \nonumber
   &=& \mathbb{E}[\int_\mathbb{R}\{\int_0^T f(t, x(t,y),u(t,y),y)\mathbb{E}[\delta_Y(y)|\mathcal{F}_t]dt +g(x(T,y),y)\mathbb{E}[\delta_Y(y)|\mathcal{F}_T]\}dy]
\end{eqnarray}
We consider the problem to find $u^{\star} \in\mathcal{A}$ such that
\begin{equation}\label{problem}
    \sup_{u\in\mathcal{A}}J(u)=J(u^{\star}).
\end{equation}

\section{A sufficient maximum principle}

The problem (\ref{problem}) is a stochastic control problem with a standard (albeit parametrized) stochastic differential equation \eqref{eq8} for the state process $x(t,y)$, but with a non-standard performance functional given by (\ref{performance}). We can solve this problem by a modified maximum principle approach, as follows:\\

Define the \emph{Hamiltonian}
 $ H:[0,T]\times\mathbb{R}\times\mathbb{R}\times\mathbb{U}\times\mathbb{R}\times\mathbb{R}\times\mathcal{R} \times \Omega \rightarrow \mathbb{R}$ by
\begin{align}\label{eq11}
&H(t,x,y,u,p,q,r)=H(t,x,y,u,p,q,r,\omega)\nonumber\\
&=\mathbb{E}[\delta_Y(y)|\mathcal{F}_t] f(t,x,u,y)+b(t,x,u,y)p + \sigma(t,x,u,y)q+\int_{\mathbb{R}}\gamma(t,x,u,y)r(y,\zeta)\nu(d\zeta).
\end{align}
Here $\mathcal{R}$ denotes the set of all functions $r(y,.) : \mathbb{R}\rightarrow  \mathbb{R}$
 such that the last integral above converges, and $p,q,r(.)$ are called the \emph{adjoint variables}.
We define the \emph{adjoint processes} $p(t,y),q(t,y),r(t,y,\zeta)$ as the solution of the $y$-parametrized BSDE
\begin{equation}\label{eq12}
    \left\{
\begin{array}{l}
    dp(t,y) = -\frac{\partial H}{\partial x}(t,y)dt +q(t,y)dB(t)+\int_{\mathbb{R}}r(t,y,\zeta)\tilde{N}(dt,d\zeta); \quad 0 \leq t\leq T\\
    p(T,y)  = g'(x(T,y),y) \mathbb{E}[\delta_Y(y)|\mathcal{F}_T]
\end{array}
    \right.
\end{equation}
Let $J(u(.,y))$ be defined by
\begin{equation}\label{J(u)2}
    J(u(.,y))=\mathbb{E}[\int_0^T f(t, x(t,y),u(t,y),y)\mathbb{E}[\delta_Y(y)|\mathcal{F}_t]dt +g(x(T,y),y)\mathbb{E}[\delta_Y(y)|\mathcal{F}_T]].
\end{equation}
Then, comparing with \eqref{performance}  we see that
\begin{equation}
J(u) = \int_{\mathbb{R}} J(u(\cdot,y)dy.
\end{equation}
Thus it suffices to maximise $J(u,y)$ over $u$ for each given parameter $y$. Hence we have transformed the original problem \eqref{problem} to the following:
\begin{problem}
 For each given $y \in \mathbb{R}$, find $u^{\star}(\cdot,y) \in\mathcal{A}$ such that
\begin{equation}\label{problem3}
    \sup_{u(.,y)\in\mathcal{A}}J(u(\cdot,y))=J(u^{\star}(\cdot,y)).
\end{equation}
\end{problem}
This is a classical (but $y$-parametrised) stochastic control problem, except  for a non-standard performance functional \eqref{J(u)2}.

To study this problem we present two maximum principles. The first is the following:
\begin{theorem}{[Sufficient maximum principle]}\\
Let $\hat{u} \in \mathcal{A}$ with associated solution $\hat{x}(t,y),\hat{p}(t,y),\hat{q}(t,y),\hat{r}(t,y,\zeta)$ of \eqref{eq8} and \eqref{eq12}. Assume that the following hold:
\begin{enumerate}
 \item $x \rightarrow g(x)$ is concave
 \item $(x,u)\rightarrow H(t,x,y,u,\widehat{p}(t,y),\widehat{q}(t,y),\hat{r}(t,y,\zeta))$ is concave for all $t,y,\zeta$
 \item $\sup_{w\in\mathbb{U}}H\big(t,\widehat{x}(t,y),w,\widehat{p}(t,y),\widehat{q}(t,y),\hat{r}(t,y,\zeta)\big)
      =H\big(t,\widehat{x}(t,y),\widehat{u}(t,y),\widehat{p}(t,y),\widehat{q}(t,y),\hat{r}(t,y,\zeta)\big)$ for all $t,y,\zeta.$
 \end{enumerate}
Then $\widehat{u}(.,y)$ is an optimal insider control for problem (\ref{problem3}).

\end{theorem}

\dproof  By considering an increasing sequence of stopping times $\tau_n$ converging to $T$, we may assume that all local integrals appearing in the computations below are martingales and have expectation 0. See \cite{OS2}. We omit the details.\\
Choose arbitrary $u(.,y)\in\mathcal{A}$, and let the corresponding solution of (\ref{eq8}) and (\ref{eq12})  be $x(t,y)$, $p(t,y)$, $q(t,y)$.
For simplicity of notation we write
$f(t,y)=f(t,x(t,y),u(t,y))$, $\widehat{f}(t,y)=f(t,\widehat{x}(t,y),\widehat{u}(t,y))$
 and similarly with $b(t,y)$, $\widehat{b}(t,y)$, $\sigma(t,y)$, $\widehat{\sigma}(t,y)$ and so on.\\
 Moreover, we write  $\widetilde{f}(t,y)=f(t,y)-\widehat{f}(t,y)$, $\widetilde{b}(t,y)=b(t,y)-\widehat{b}(t,y)$, $\widetilde{x}(t,y)=x(t,y)-\widehat{x}(t,y)$.\\
  Consider
 \begin{equation*}
    J(u(.,y))-J(\widehat{u}(.,y))=I_1+ I_2,
 \end{equation*}
 where
 \begin{equation}\label{I_1I_2}
    I_1=\mathbb{E}[\int_0^T\{f(t,y)-\widehat{f}(t,y)\}\mathbb{E}[\delta_Y(y)|\mathcal{F}_t]dt], \quad I_2=\mathbb{E}[\{g(x(T,y))-g(\widehat{x}(T,y))\}\mathbb{E}[\delta_Y(y)|\mathcal{F}_T]].
 \end{equation}
  By the definition of $H$ we have
  \begin{eqnarray}\label{II1}
    I_1 &=& \mathbb{E}[\int_0^T\{H(t,y)-\widehat{H}(t,y)-\widehat{p}(t,y)\widetilde{b}(t,y) - \widehat{q}(t,y)\widetilde{\sigma}(t,y)\\ \nonumber
    &-&\int_{\mathbb{R}}\hat{r}(t,y,\zeta)\tilde{\gamma}(t,y,\zeta)\nu(d\zeta)\}dt].
  \end{eqnarray}

   Since $g$ is concave we have by (\ref{eq12})
 \begin{align}\label{II_2}
   I_2 &\leq\mathbb{E}[g'(\widehat{x}(T,y))\mathbb{E}[\delta_Y(y)|\mathcal{F}_t](x(T,y)-\widehat{x}(T,y))]=\mathbb{E}[\widehat{p}(T,y)\widetilde{x}(T,y)] \\ \nonumber
    &= \mathbb{E}[\int_0^T \widehat{p}(t,y) d\widetilde{x}(t,y)+\int_0^T\widetilde{x}(t,y)d\widehat{p}(t,y)+\int_0^Td[\hat{p} , \tilde{x}]_t] \\ \nonumber
    &=\mathbb{E}[\int_0^T \widehat{p}(t,y) (\widetilde{b}(t,y)dt+\widetilde{\sigma}(t,y)dB(t)+\int_{\mathbb{R}}\tilde{\gamma}(t,y,\zeta)\tilde{N}(dt,d\zeta)) \\ \nonumber
    &- \int_0^T\frac{\partial \widehat{H}}{\partial x}(t,y)\widetilde{x}(t,y)dt+\int_0^T\widehat{q}(t,y)\widetilde{x}(t,y)dB(t)+\int_0^T\int_{\mathbb{R}}\widetilde{x}(t,y)\hat{r}(t,y,\zeta)\tilde{N}(dt,d\zeta)\\ \nonumber
    &+\int_0^T\widetilde{\sigma}(t,y)\widehat{q}(t,y)dt+\int_0^T\int_{\mathbb{R}}\tilde{\gamma}(t,y,\zeta)\hat{r}(t,y,\zeta)\nu(d\zeta)dt
    +\int_0^T\int_{\mathbb{R}}\tilde{\gamma}(t,y,\zeta)\hat{r}(t,y,\zeta)\tilde{N}(dt,d\zeta)]\\  \nonumber
    &= \mathbb{E}[\int_0^T \widehat{p}(t,y) \widetilde{b}(t,y)dt -\int_0^T\frac{\partial \widehat{H}}{\partial x}(t,y)\widetilde{x}(t,y)dt
    +\int_0^T\widetilde{\sigma}(t)\widehat{q}(t)dt+\int_0^T\int_{\mathbb{R}}\tilde{\gamma}(t,y,\zeta)\hat{r}(t,y,\zeta)\nu(d\zeta)dt].
     \end{align}

Adding (\ref{II1}) - (\ref{II_2}) we get, by concavity of $H$,
\begin{eqnarray*}
  J(u(.,y))-J(\widehat{u}(.,y)) &\leq&\mathbb{E}[\int_0^T \{H(t,y)-\widehat{H}(t,y)-\frac{\partial \hat{H}}{\partial x}(t,y)\widetilde{x}(t,y)\}dt ] \\ \nonumber
   &\leq&  \mathbb{E}[\int_0^T \frac{\partial \hat{H}}{\partial u}(t,y)\widetilde{u}(t,y) dt] \\ \nonumber
   &\leq& 0,
\end{eqnarray*}
since $u(.,y)=\widehat{u}(.,y)$ maximizes $\widehat{H}(.,y)$ at $t$.
\fproof

\section{A necessary maximum principle}
We proceed to establish a corresponding necessary maximum principle. For this, we do not need concavity conditions, but in stead we need the following assumptions about the set of admissible control values:\\
\begin{itemize}
\item
$A_1$. For all $t_0\in [0,T]$ and all bounded $\mathcal{H}_{t_0}$-measurable random variables $\alpha(y,\omega)$, the control
$\theta(t,y, \omega) := \mathbf{1}_{[t_0,T ]}(t)\alpha(y,\omega)$ belongs to $\mathcal{A}$.\\
\item
$A_2$. For all $u; \beta_0 \in\mathcal{A}$ with $\beta_0(t,y) \leq K < \infty$ for all $t,y$  define
\begin{equation}\label{delta}
    \delta(t,y)=\frac{1}{2K}dist((u(t,y),\partial\mathbb{U})\wedge1 > 0
\end{equation}
and put
\begin{equation}\label{beta(t,y)}
    \beta(t,y)=\delta(t,y)\beta_0(t,y).
\end{equation}
Then the control
\begin{equation*}
    \widetilde{u}(t,y)=u(t,y) + a\beta(t,y) ; \quad t \in [0,T]
\end{equation*}
belongs to $\mathcal{A}$ for all $a \in (-1, 1)$.\\
\item
$A3$. For all $\beta$ as in (\ref{beta(t,y)}) the derivative process
\begin{equation*}
    \chi(t,y):=\frac{d}{da}x^{u+a\beta}(t,y)|_{a=0}
\end{equation*}
exists, and belong to $\mathbf{L}^2(\lambda\times \mathbf{P})$ and
\begin{equation}\label{d chi}
    \left\{
\begin{array}{l}
    d\chi(t,y) = [\frac{\partial b}{\partial x}(t,y)\chi(t,y)+\frac{\partial b}{\partial u}(t,y)\beta(t,y)]dt+[\frac{\partial \sigma}{\partial x}(t,y)\chi(t,y)+\frac{\partial\sigma}{\partial u}(t,y)\beta(t,y)]d B(t) \\
    +\int_{\mathbb{R}}[\frac{\partial \gamma}{\partial x}(t,y,\zeta)\chi(t,y)+\frac{\partial \gamma}{\partial u}(t,y,\zeta)\beta(t,y)]\tilde{N}(dt,d\zeta)\\
    \chi(0,y)  = \frac{d}{da}x^{u+a\beta}(0,y)|_{a=0} = 0.
   \end{array}
    \right.
\end{equation}
\end{itemize}

\begin{theorem}{[Necessary maximum principle]} \\
Let $\hat{u} \in \mathcal{A}$. Then the following are equivalent:
\begin{enumerate}
\item $\frac{d}{da}J((\hat{u}+a\beta)(.,y))|_{a=0}=0$ for all bounded $\beta \in \mathcal{A}$ of the form (\ref{beta(t,y)}).
\item $\frac{\partial H}{\partial u}(t,y)_{u=\hat{u}}=0$ for all $t\in[0,T].$
\end{enumerate}
\end{theorem}

\dproof
For simplicity of notation we write $u$ instead of $\hat{u}$ in the following. \\
By considering an increasing sequence of stopping times $\tau_n$ converging to $T$, we may assume that all local integrals appearing in the computations below are martingales and have expectation 0. See \cite{OS2}. We omit the details.\\
We can write
$$\frac{d}{da}J((u+a\beta)(.,y))|_{a=0}=I_1+I_2$$\\
where
$$I_1=\frac{d}{da}\mathbb{E}[\int_0^Tf(t,x^{u+a\beta}(t,y),u(t,y)+a\beta(t,y),y)\mathbb{E}[\delta_Y(y)|\mathcal{F}_t]dt]|_{a=0}$$\\
and
$$I_2=\frac{d}{da}\mathbb{E}[g(x^{u+a\beta}(T,y),y)\mathbb{E}[\delta_Y(y)|\mathcal{F}_T]]|_{a=0}.$$
By our assumptions on $f$ and $g$ and by \eqref{eq12} we have
\begin{equation}\label{iii1}
    I_1=\mathbb{E}[\int_0^T\{\frac{\partial f}{\partial x}(t,y)\chi(t,y)+\frac{\partial f}{\partial u}(t,y)\beta(t,y)\}\mathbb{E}[\delta_Y(y)|\mathcal{F}_t]dt]
\end{equation}
\begin{equation}\label{iii2}
    I_2=\mathbb{E}[g'(x(T,y),y)\chi(T,y)\mathbb{E}[\delta_Y(y)|\mathcal{F}_T]]=\mathbb{E}[p(T,y)\chi(T,y)]
\end{equation}
By the It\^{o} formula
\begin{eqnarray}\label{iii22}
   I_2&=& \mathbb{E}[p(T,y)\chi(T,y)]=\mathbb{E}[\int_0^Tp(t,y)d\chi(t,y)+\int_0^T\chi(t,y)dp(t,y)+\int_0^Td[\chi,p](t,y)] \\ \nonumber
   &=& \mathbb{E}[\int_0^Tp(t,y)\{\frac{\partial b}{\partial x}(t,y)\chi(t,y)+\frac{\partial b}{\partial u}(t,y)\beta(t,y)\}dt\\ \nonumber
   &+&\int_0^Tp(t,y)\{\frac{\partial \sigma}{\partial x}(t,y)\chi(t,y)+\frac{\partial\sigma}{\partial u}(t,y)\beta(t,y)\}dB(t) \\ \nonumber
   &+&\int_0^T\int_{\mathbb{R}}p(t,y)\{\frac{\partial \gamma}{\partial x}(t,y,\zeta)\chi(t,y)+\frac{\partial \gamma}{\partial u}(t,y,\zeta)\beta(t,y)\}\tilde{N}(dt,d\zeta)\\ \nonumber
   &-&\int_0^T\chi(t,y)\frac{\partial H}{\partial x}(t,y) dt+\int_0^T\chi(t,y)q(t,y)dB(t)+\int_0^T\int_{\mathbb{R}}\chi(t,y)r(t,y,\zeta) \tilde{N}(dt,d\zeta)  \\ \nonumber
     &+&\int_0^Tq(t,y) \{\frac{\partial \sigma}{\partial x}(t,y)\chi(t,y)+\frac{\partial\sigma}{\partial u}(t,y)\beta(t,y)\}dt \\ \nonumber
   &+&\int_0^T\int_{\mathbb{R}}\{\frac{\partial \gamma}{\partial x}(t,y,\zeta)\chi(t,y)+\frac{\partial \gamma}{\partial u}(t,y,\zeta)\beta(t,y)\}r(t,y,\zeta)\nu(\zeta)dt]\\ \nonumber
  &=& \mathbb{E}[\int_0^T\chi(t,y)\{p(t,y)\frac{\partial b}{\partial x}(t,y)+q(t,y)\frac{\partial \sigma}{\partial x}(t,y)-\frac{\partial H}{\partial x}(t,y)+\int_{\mathbb{R}}\frac{\partial \gamma}{\partial x}(t,y,\zeta)r(t,y,\zeta)\nu(d\zeta)\}dt  \\ \nonumber
   &+& \int_0^T\beta(t,y)\{p(t,y)\frac{\partial b}{\partial u}(t,y)+q(t,y)\frac{\partial\sigma}{\partial u}(t,y)+\int_{\mathbb{R}}\frac{\partial \gamma}{\partial u}(t,y,\zeta)r(t,y,\zeta)\nu(d\zeta)\}dt]\\ \nonumber
&=&\mathbb{E}[-\int_0^T\chi(t,y)\frac{\partial f}{\partial x}\mathbb{E}[\delta_Y(y)|\mathcal{F}_t]dt+\int_0^T\{\frac{\partial H}{\partial u}(t,y)-\frac{\partial f}{\partial u}(t,y)\mathbb{E}[\delta_Y(y)|\mathcal{F}_t]\}\beta(t,y)dt]\\ \nonumber
&=& -I_1+\mathbb{E}[\int_0^T\frac{\partial H}{\partial u}(t,y)\beta(t,y)dt].
\end{eqnarray}

Summing (\ref{iii1}) and (\ref{iii22}) we get
\begin{equation*}
    \frac{d}{da}J((u+a\beta)(.,y))|_{a=0}=I_1+I_2=\mathbb{E}[\int_0^T\frac{\partial H}{\partial u}(t,y)\beta(t,y)dt].
\end{equation*}
we conclude that
\begin{equation*}
    \frac{d}{da}J((u+a\beta)(.,y))|_{a=0}=0
\end{equation*}
if and only if $\mathbb{E}[\int_0^T\frac{\partial H}{\partial u}(t,y)\beta(t,y)dt]=0$ for all bounded $\beta\in\mathcal{A}$ of the form (\ref{beta(t,y)}).\\

\noindent In particular, applying this to $\beta(t,y) = \theta(t,y)$ as in $A1$, we get that this is again equivalent to
\begin{equation*}
    \frac{\partial H}{\partial u}(t,y)=0  \text{ for all } t\in[0,T].
\end{equation*}
\fproof

\section{Applications}
In the following we assume that

\begin{equation}
\mathbb{E}[\int_0^T \{\mathbb{E}[D_t \delta_Y(y)|\mathcal{F}_t]^2 +\int_{\mathbb{R}}\mathbb{E}[D_{t,z}\delta_Y(y)|\mathcal{F}_t]^2 \nu(dz) \} dt] < \infty.
\end{equation}

\subsection{Utility maximization for an insider, part 1 (N=0)}
Consider a financial market where the unit price $S_0(t)$ of the risk free asset is
\begin{equation}\label{riskfree}
    S_0(t)=1, \quad t\in[0, T]
\end{equation}
and the unit price process $S(t)$ of the risky asset has no jumps and is given by

\begin{align}\label{eq5.2}
\begin{cases}
 dS(t) &= S(t) [b_0(t,Y) dt + \sigma_0(t,Y) dB(t)]
 ; \quad t\in[0, T]\\
    S(0)&>0.
\end{cases}
\end{align}

Then the wealth process $X(t)=X^{\Pi}(t)$ associated to a portfolio $u(t)=\Pi(t)$, interpreted as the fraction of the wealth invested in the risky asset at time $t$, is given by
\begin{align}\label{eq5.3}
  \begin{cases}
dX(t)&=\Pi(t)X(t)[b_0(t,Y)dt+\sigma_0(t,Y)dB(t)]
; \quad t\in[0, T]\\
X(0)&=x_0>0.
\end{cases}
  \end{align}

Let $U$ be a given utility function. We want to find $\Pi^{\ast}\in \mathcal{A}$ such that
\begin{equation}\label{eq17}
J(\Pi^*) =\sup_{\Pi\in\mathcal{A}} J(\Pi),
\end{equation}
where
\begin{equation}\label{eq18}
J(\Pi):=  \mathbb{E}[U(X^\Pi(T))].
\end{equation}
Note that, in terms of our process $x(t,y)$ we have
\begin{align}\label{Wealth}
\begin{cases}
dx(t,y)=\pi(t,y)x(t,y)[b_0(t,y)dt+\sigma_0(t,y)dB(t)]
; \quad t\in[0, T]\\
x(0,y)=x_0>0
\end{cases}
\end{align}
and the performance functional gets the form
\begin{equation*}
J(\pi)= \mathbb{E}[U(x(T,y)) \mathbb{E}[\delta_Y(y)|\mathcal{F}_T]].
\end{equation*}
This is a problem of the type investigated in the previous sections (in the special case with no jumps) and we can apply the results there to solve it, as follows:\\
The Hamiltonian gets the form, with $u=\pi$,
\begin{equation}\label{eq19}
H(t,x,y,\pi,p,q)= \pi x [b_0(t,y)p + \sigma_0(t,y)q]
\end{equation}
while the BSDE for the adjoint processes becomes
\begin{align}\label{eq20}
\begin{cases}
dp(t,y) &= - \pi(t,y)  [b_0(t,y)p(t,y) + \sigma_0(t,y)q(t,y)]dt
+ q(t,y) dB(t)
; \quad t\in[0, T]\\
p(T,y)&= U'(x(T,y))\mathbb{E}[\delta_Y(y)|\mathcal{F}_T]
\end{cases}
\end{align}

Since the Hamiltonian $H$ is a linear function of $\pi$, it can have a finite maximum over all $\pi$ only if
\begin{equation}\label{eq21}
 x(t,y)  [b_0(t,y)p(t,y) + \sigma_0(t,y)q(t,y)]=0
\end{equation}

Substituted into \eqref{eq20} this gives
\begin{align}\label{eq22}
\begin{cases}
dp(t,y) &= q(t,y) dB(t)
 \\
p(T,y)&= U'(x(T))\mathbb{E}[\delta_Y(y)|\mathcal{F}_T]
\end{cases}
\end{align}

If we assume that, for all $t,y$,
\begin{equation}\label{eq23}
 x(t,y) > 0,
\end{equation}
\noindent then we get from \eqref{eq21} that
\begin{equation}
q(t,y) = -  \frac{b_0(t,y)}{\sigma_0(t,y)} p(t,y).
\end{equation}

\ Substituting this into \eqref{eq22}, we get the equation
 \begin{align}\label{eq25}
\begin{cases}
dp(t,y) &= - \frac{b_0(t,y)}{\sigma_0(t,y)} p(t,y) dB(t)\\
p(T,y)&= U'(x(T,y))\mathbb{E}[\delta_Y(y)|\mathcal{F}_T]
\end{cases}
\end{align}

Thus we obtain that
\begin{equation}\label{eq26}
p(t,y)=p(0,y) \exp(-\int_0^t \frac{b_0(s,y)}{\sigma_0(s,y)}dB(s) - \frac{1}{2} \int_0^t (\frac{b_0(s,y)}{\sigma_0(s,y)})^{2}ds),
\end{equation}
for some, not yet determined, constant $p(0,y)$.
In particular, if we put $t=T$ and use \eqref{eq25} we get
\begin{equation}\label{eq26b}
U'(x(T,y))\mathbb{E}[\delta_Y(y)|\mathcal{F}_T]=p(0,y) \exp(-\int_0^T \frac{b_0(s,y)}{\sigma_0(s,y)}dB(s) - \frac{1}{2} \int_0^T (\frac{b_0(s,y)}{\sigma_0(s,y)})^{2}ds).
\end{equation}

\noindent To make this more explicit, we proceed as follows:\\

Define
\begin{equation}
M(t,y) := \mathbb{E}[\delta_Y(y)|\mathcal{F}_t]
\end{equation}
Then by the generalized Clark-Ocone theorem
\begin{equation}
\begin{cases}
dM(t,y) = \mathbb{E}[D_t \delta_Y(y) | \mathcal{F}_t]dB(t)= \Phi(t,y)M(t,y) dB(t)\\
M(0,y)=1
\end{cases}
\end{equation}
where
\begin{equation}\label{eq18a}
\Phi(t,y)= \frac{\mathbb{E}[D_t\delta_Y(y)|\mathcal{F}_t]}{\mathbb{E}[\delta_Y(y)|\mathcal{F}_t]}
\end{equation}
Solving this SDE for $M(t)$ we get
\begin{equation}
M(t) = \exp( \int_0^t \Phi(s,y)dB(s) -\frac{1}{2} \int_0^t \Phi^2(s,y) ds).
\end{equation}
Substituting this into \eqref{eq26b} we get
\begin{align}\label{eq20a}
&U'(x(T,y))= p(0,y) \exp \big( - \int_0^T \{\Phi(s,y) +\frac{b_0(s,y)}{\sigma_0(s,y)} \}dB(s)\nonumber\\
&+\frac{1}{2} \int_0^T\{\Phi^2(s,y) -\frac{b_0^2(s,y)}{\sigma_0^2(s,y)}\} ds \big)=:p(0,y)\Gamma(T,y) .
\end{align}
i.e.,
\begin{equation}\label{eq34}
    x(T,y)=I(c\Gamma(T,y))
\end{equation}
where
\begin{equation}
I=(U')^{-1} \text{  and } c=p(0,y).
\end{equation}
It remains to find $c$. We can write the differential stochastic equation of $x(t,y)$  as

\begin{equation}\label{eq35}
    \left\{
  \begin{array}{lr}
dx(t,y)=\pi(t,y)x(t,y)[b_0(t,y)dt+\sigma_0(t,y)dB(t)]\\
x(T,y)=I(c\Gamma(T,y))
  \end{array}
    \right.
\end{equation}

If we define
\begin{equation}\label{eq36}
    z(t,y)= \pi(t,y)x(t,y)\sigma_0(t,y)
\end{equation}
then equation (\ref{eq35}) becomes the linear BSDE
\begin{equation}\label{eq37}
    \left\{
  \begin{array}{lr}
dx(t,y)=\frac{z(t,y)b_0(t,y)}{\sigma_0(t,y)}dt+z(t,y)dB(t)\\
x(T,y)=I(c\Gamma(T,y))
  \end{array}
    \right.
\end{equation}
in the unknown $(x(t,y),z(t,y))$.
The solution of this BSDE is
\begin{equation}\label{eq38}
    x(t,y)=\frac{1}{\Gamma_0(t,y)}\mathbb{E}[I(c\Gamma(T,y))\Gamma_0(T,y)|\mathcal{F}_t],
\end{equation}
where
\begin{equation}
\Gamma_0(t,y)=\exp\{-\int_0^t\frac{b_0(s,y)}{\sigma_0(s,y)}dB(s)-\frac{1}{2}\int_0^t(\frac{b_0(s,y)}{\sigma_0(s,y)})^2ds\}.
\end{equation}
In particular,
\begin{equation}\label{eq39}
    x_0=x(0,y)=\mathbb{E}[I(c\Gamma(T,y))\Gamma_0(T,y)].
\end{equation}
This is an equation which (implicitly) determines the value of $c$.
When $c$ is found, we have the optimal terminal wealth $x(T,y)$ given by (\ref{eq37}). Solving the
resulting BSDE for  $z(t,y)$, we get the corresponding optimal portfolio $\pi(t,y)$ by (\ref{eq36}).
We summarize what we have proved in the following theorem:

\begin{theorem}
The optimal portfolio $\Pi^*(t)$ for the insider portfolio problem \eqref{eq17} is given by
\begin{equation}
\Pi^*(t) = \int_{\mathbb{R}}\pi^*(t,y) \delta_Y(y) dy=\pi^*(t,Y),
\end{equation}
where
\begin{equation}
\pi^*(t,y)=\frac{z(t,y)}{x(t,y)\sigma_0(t,y)}
\end{equation}
with $x(t,y), z(t,y)$ given as the solution of the BSDE \eqref{eq37} and $c=p(0,y)$ given by \eqref{eq39}.
\end{theorem}

\subsection{The logarithmic utility case (N=0)}
We now look at the special case when $U$ is the \emph{logarithmic utility}, i.e.,
\begin{equation}
U(x) = \ln x; \text{  } x > 0.
\end{equation}
Recall the equation for $x(t,y)$:
\begin{align}\label{Wealth}
\begin{cases}
dx(t,y)=\pi(t,y)x(t,y)[b_0(t,y)dt+\sigma_0(t,y)dB(t)]\\
x(0,y)=x_0>0
\end{cases}
\end{align}
By the It\^{o} formula for forward integrals, we get that the solution
of this equation is
\begin{equation}\label{solution x(t,y)}
    x(t,y) = x_0 \exp\{\int_0^t[\pi(s,y)b_0(s,y)-\frac{1}{2}\pi^2(s,y)\sigma_0^2(s,y)]ds + \int_0^t\pi(s,y)\sigma_0(s,y)dB(s)\}.
\end{equation}
Therefore,
\begin{align}\label{U'(x(T,y)}
& U' (x(T,y)) \nonumber\\
&=\frac{1}{x(T,y)}= \frac{1}{x}\exp\{-\int_0^T[\pi(s,y)b_0(s,y)-\frac{1}{2}\pi^2(s,y)\sigma_0^2(s,y)]ds - \int_0^T\pi(s,y)\sigma_0(s,y)dB(s)\}  \nonumber\\
\end{align}
Comparing with \eqref{eq20a} we try to choose $\pi(s,y)$ such that
\begin{align}\label{eq41a}
& \frac{1}{x}\exp\{- \int_0^T\pi(s,y)\sigma_0(s,y)dB(s)-\int_0^T[\pi(s,y)b_0(s,y)-\frac{1}{2}\pi^2(s,y)\sigma_0^2(s,y)]ds \}\nonumber\\
&= p(0,y) \exp\{- \int_0^T \{\Phi(s,y) +\frac{b_0(s,y}{\sigma_0(s,y)} \}dB(s)
+\frac{1}{2} \int_0^T\int_0^t \{\Phi^2(s,y) -\frac{b_0^2(s,y)}{\sigma_0^2(s,y)}\} ds\}
\end{align}
Thus we try to put
\begin{equation}
p(0,y)= \frac{1}{x}
\end{equation}
and choose $\pi(s,y)$ such that, using \eqref{eq18a},
\begin{align}
\pi(s,y)\sigma_0(s,y) = \Phi(s,y) +\frac{b_0(s,y)}{\sigma_0(s,y)}
=\frac {\mathbb{E}[D_s\delta_Y(y)|\mathcal{F}_s]}{\mathbb{E}[\delta_Y(y)|\mathcal{F}_s]}+\frac{b_0(s,y)}{\sigma_0(s,y)}
\end{align}
This gives
\begin{equation}\label{pi(s,y)}
    \pi(s,y)=\frac{b_0(s,y)}{\sigma^2_0(s,y)}+\frac {\mathbb{E}[D_s\delta_Y(y)|\mathcal{F}_s]}{\sigma_0(s,y)\mathbb{E}[\delta_Y(y)|\mathcal{F}_s]}
\end{equation}
We now verify that with this choice of $\pi(s,y)$, also the other two terms on the exponents on each side of \eqref{eq41a} coincide, i.e. that
\begin{equation}
\pi(s,y)b_0(s,y)-\frac{1}{2}\pi^2(s,y)\sigma_0^2(s,y)=\Phi^2(s,y) -\frac{b_0^2(s,y)}{\sigma_0^2(s,y)}
\end{equation}

Thus we have proved the following result, which has been obtained earlier in \cite{OR1} by a different method:

\begin{theorem}
The optimal portfolio $\Pi = \Pi^*$ with respect to logarithmic utility for an insider in the market (\ref{riskfree})-\eqref{eq5.2} and with the inside information (\ref{H_t}) is given by
\begin{equation}\label{pi(s,y)}
    \Pi^*(s)=\frac{b_0(s,Y)}{\sigma^2_0(s,Y)}+\frac {\mathbb{E}[D_s\delta_Y(y)|\mathcal{F}_s]_{y=Y}}{\sigma_0(s,Y)\mathbb{E}[\delta_Y(y)|\mathcal{F}_s]}_{y=Y}, \quad 0 \leq s \leq T < T_0.
\end{equation}
\end{theorem}

Substituting (\ref{eq5.50}) and (\ref{eq5.51}) in (\ref{pi(s,y)}) we obtain:
\begin{coro}
Suppose that Y is Gaussian of the form (\ref{eq5.47}). Then the optimal insider
portfolio is given by
\begin{equation}\label{eq5.52}
    \Pi^*(s)=\frac{b_0(s,Y(T_0))}{\sigma^2_0(s,Y(T_0))}+\frac {(Y(T_0)-Y(s))\beta(s)}{\sigma_0(s,Y(T_0))\|\beta\|^2_{[s,T_0]}}, \quad 0 \leq s \leq T < T_0.
\end{equation}
\end{coro}
In particular, if $Y = B(T_0)$ we get the following result, which was also proved in \cite{PK}, in the case when the coefficients do not depend on $Y$:
\begin{coro}
Suppose that $Y = B(T_0)$. Then the optimal insider portfolio is given by
\begin{equation}\label{eq5.52}
    \Pi^*(s)=\frac{b_0(s,B(T_0))}{\sigma^2_0(s,B(T_0))}+\frac {B(T_0)-B(s)}{\sigma_0(t,B(T_0))(T_0-s)},\quad 0 \leq s \leq T < T_0.
\end{equation}
\end{coro}

\subsection{Utility maximization for an insider, part 2 (Poisson process case)}
Let $N(t)$ be the Poisson process with intensity $\lambda > 0$. Consider a financial market where the unit price $S_0(t)$ of the risk free asset is
\begin{equation}\label{riskfree'}
    S_0(t)=1, \quad t\in[0, T]
\end{equation}
and the unit price process $S(t)$ of the risky asset has no jumps and is given by

\begin{align}\label{eq5.2'}
\begin{cases}
 dS(t) &= S(t) [b_0(t,Y) dt + \gamma_0(t,Y) d\tilde{N}(t)]
 ; \quad t\in[0, T_0]\\
    S(0)&>0.
\end{cases}
\end{align}
where $\tilde{N}=N(t) -\lambda t$ is the compensated Poisson process with parameter $\lambda$.
In this case the L\'  evy measure is $\nu(d\zeta)=\lambda\delta_1(d\zeta)$ since the jumps are with size $1$.
As before $Y$ is a given random variable whose value is known to the trader at any time $t \leq T < T_0$.
Then the wealth process $X(t)=X^{\Pi}(t)$ associated to a portfolio $u(t)=\Pi(t)$, interpreted as the fraction of the wealth invested in the risky asset at time $t$, is given by
\begin{align}\label{eq5.3'}
  \begin{cases}
dX(t)&=\Pi(t)X(t)[b_0(t,Y)dt+\gamma_0(t,Y)d\tilde{N}(t)]
; \quad t\in[0, T]\\
X(0)&=x_0>0.
\end{cases}
  \end{align}

Let $U$ be a given utility function. We want to find $\Pi^{\ast}\in \mathcal{A}$ such that
\begin{equation}\label{eq17'}
J(\Pi^*) =\sup_{\Pi\in\mathcal{A}} J(\Pi),
\end{equation}
where
\begin{equation}\label{eq18'}
J(\Pi):=  \mathbb{E}[U(X^\Pi(T))].
\end{equation}
Note that, in terms of our process $x(t,y)$ we have
\begin{align}\label{Wealth}
\begin{cases}
dx(t,y)=\pi(t,y)x(t,y)[b_0(t,y)dt+\gamma_0(t,y)d\tilde{N}(t)]
; \quad t\in[0, T]\\
x(0,y)=x_0(y)>0,
\end{cases}
\end{align}
which has the solution
\begin{align}\label{eqWealth2}
x(t,y)& = x_0(y) \exp \big( \int_0^t \{\pi(s,y) b_0(s,y) + \lambda \ln(1+\pi(s,y)\gamma_0(s,y)) - \lambda \pi(s,y) \gamma_0(s,y) \} ds \nonumber\\
&+ \int_0^t  \ln(1+ \pi(s,y) \gamma_0(s,y)) d\tilde{N}(s) \big).
\end{align}
The performance functional gets the form
\begin{equation}
J(\pi)=\mathbb{E}[ U(x(T,y)) \mathbb{E}[\delta_Y(y)|\mathcal{F}_T]],
\end{equation}
where
\begin{equation}
\Pi(t)=\pi(t,Y).
\end{equation}

This is a problem of the type investigated in the previous sections, in the special case when all the jumps are of size $1$. Then the L\' evy measure gets the form $\nu(d\zeta)=\lambda\delta_1$ and when we apply the results from the previous section we get:\\
The Hamiltonian becomes, with $u=\pi$,
\begin{equation}\label{eq19'}
H(t,x,y,\pi,p,r)= \pi x [b_0(t,y)p +\lambda \gamma_0(t,y)r(y,1)]
\end{equation}
while the BSDE for the adjoint processes becomes
\begin{align}\label{eq20'}
\begin{cases}
dp(t,y) &= - \pi(t,y)  [b_0(t,y)p(t,y) + \lambda \gamma_0(t,y) r(t,y,1) ]dt
+ r(t,y,1) d\tilde{N}(t)
; \quad t\in[0, T]\\
p(T,y)&= U'(x(T,y))\mathbb{E}[\delta_Y(y)|\mathcal{F}_T]
\end{cases}
\end{align}

Since the Hamiltonian $H$ is a linear function of $\pi$, it can have a finite maximum over all $\pi$ only if
\begin{equation}\label{eq21'}
 x(t,y)  [b_0(t,y)p(t,y) + \gamma_0(t,y)\lambda r(t,y,1)]=0
\end{equation}

Substituted into \eqref{eq20'} this gives
\begin{align}\label{eq22'}
\begin{cases}
dp(t,y) &= r(t,y,1)\tilde{N}(dt) \\
p(T,y)&= U'(x(T))\mathbb{E}[\delta_Y(y)|\mathcal{F}_T]
\end{cases}
\end{align}

If we assume that, for all $t,y$,
\begin{equation}\label{eq23'}
 x(t,y) > 0,
\end{equation}
\noindent then we get from \eqref{eq21'} that
\begin{equation}
r(t,y,1) = -  \frac{b_0(t,y)}{\lambda\gamma_0(t,y)} p(t,y).
\end{equation}

\ Substituting this into \eqref{eq22'}, we get the equation
 \begin{align}\label{eq25'}
\begin{cases}
dp(t,y) &= - \frac{b_0(t,y)}{\lambda \gamma_0(t,y)} p(t,y) d\tilde{N}(t)\\
p(T,y)&= U'(x(T,y))\mathbb{E}[\delta_Y(y)|\mathcal{F}_T]
\end{cases}
\end{align}

Thus we obtain that (see e.g. \cite{OS1},Example 1.15)
\begin{equation}\label{eq26'}
p(t,y)=p(0,y) \exp(\int_0^t \ln[1-\frac{b_0(s,y)}{\lambda\gamma_0(s,y)}]d\tilde{N}(s) + \lambda \int_0^t (\ln[1-\frac{b_0(s,y)}{\lambda\gamma_0(s,y)}]+\frac{b_0(s,y)}{\lambda\gamma_0(s,y)})ds),
\end{equation}
for some, not yet determined, constant $p(0,y)$.
In particular, if we put $t=T$ and use \eqref{eq25'} we get
\begin{align}\label{eq26b'}
U'(x(T,y))\mathbb{E}[\delta_Y(y)|\mathcal{F}_T]&=p(0,y) \exp \big( \int_0^T \ln[1-\frac{b_0(s,y)}{\lambda\gamma_0(s,y)}]d\tilde{N}(s) \nonumber\\
&+ \lambda \int_0^T (\ln[1-\frac{b_0(s,y)}{\lambda\gamma_0(s,y)}]+\frac{b_0(s,y)}{\lambda\gamma_0(s,y)})ds \big).
\end{align}

\noindent To make this more explicit, we proceed as follows:\\

Define

\begin{equation}
M(t,y) := \mathbb{E}[\delta_Y(y)|\mathcal{F}_t]
\end{equation}
Then by the generalized Clark-Ocone theorem
\begin{equation}
\begin{cases}
dM(t,y) = \mathbb{E}[D_{t,1} \delta_Y(y) | \mathcal{F}_t]d\tilde{N}(t)= \Psi(t,y)M(t,y) d\tilde{N}(t)\\
M(0,y)=1
\end{cases}
\end{equation}
where
\begin{equation}\label{eq18a'}
\Psi(t,y)= \frac{\mathbb{E}[D_{t,1}\delta_Y(y)|\mathcal{F}_t]}{\mathbb{E}[\delta_Y(y)|\mathcal{F}_t]}
\end{equation}
Solving this SDE for $M(t)$ we get
\begin{equation}
M(t,y) = \exp( \int_0^t \ln(1+\Psi(s,y))d\tilde{N}(s) +\lambda\int_0^t [\ln(1+\Psi(s,y))-\Psi(s,y)] ds).
\end{equation}
Substituting this into \eqref{eq26b'} we get
\begin{align}\label{eq20a'}
U'(x(T,y))=&p(0,y) \exp \big(  \int_0^T [\ln(1-\frac{b_0(s,y)}{\lambda\gamma_0(s,y)})-\ln(1+\Psi(s,y))]d\tilde{N}(s)\nonumber\\
&+\lambda\int_0^T\{[\ln(1-\frac{b_0(s,y)}{\lambda\gamma_0(s,y)})+\frac{b_0(s,y)}{\lambda\gamma_0(s,y)}]-[\ln(1+\Psi(s,y))\nonumber\\
&-\Psi(s,y)]\} ds \big)=:p(0,y)\Gamma(T,y),
\end{align}\\
i.e.,
\begin{equation}
 x(T,y)=I(c\Gamma(T,y))
\end{equation}
where
\begin{equation}
I=(U')^{-1} \text{ and } c=p(0,y).
\end{equation}
It remains to find $c$.
We can write the differential stochastic equation of $x(t,y)$  as
\begin{equation}\label{eq35'}
    \left\{
  \begin{array}{lr}
dx(t,y)=\pi(t,y)x(t,y)[b_0(t,y)dt+\gamma_0(t,y)d\tilde{N}(t)]\\
x(T,y)=I(c\Gamma(T,y))
  \end{array}
    \right.
\end{equation}

If we define
\begin{equation}\label{eq36'}
    k(t,y)):= \pi(t,y)x(t,y)\gamma_0(t,y)
\end{equation}
then equation (\ref{eq35'}) becomes the BSDE
\begin{equation}\label{eq37'}
    \left\{
  \begin{array}{lr}
dx(t,y)=\frac{b_0(t,y)}{\lambda\gamma_0(t,y)}k(t,y)\lambda dt+k(t,y)d\tilde{N}(t)\\
x(T,y)=I(c\Gamma(T,y))
  \end{array}
    \right.
\end{equation}
in the unknown $(x(t,y),k(t,y))$.
The solution of this BSDE is
\begin{equation}\label{eq38'}
    x(t,y)=\frac{1}{\Gamma_0(t,y)}\mathbb{E}[I(c\Gamma(T,y))\Gamma_0(T,y)|\mathcal{F}_t],
\end{equation}
where
\begin{equation}
\Gamma_0(t,y)=\exp\{\int_0^t[\ln(1-\frac{b_0(s,y)}{\lambda\gamma_0(s,y)})+\frac{b_0(s,y)}{\lambda\gamma_0(s,y)}]\lambda ds+\int_0^t\ln(1-\frac{b_0(s,y)}{\lambda\gamma_0(s,y)})\tilde{N}(ds)\}.
\end{equation}
In particular,
\begin{equation}\label{eq39'}
    x_0=x(0,y)=\mathbb{E}[I(c\Gamma(T,y))\Gamma_0(T,y)].
\end{equation}

This is an equation which (implicitly) determines the value of $c$.
When $c$ is found, we have the optimal terminal wealth $x(T,y)$ given by (\ref{eq37'}). Solving the
resulting BSDE for  $k(t,y)$, we get the corresponding optimal portfolio $\pi(t,y)$ by (\ref{eq36'}).
We summarize what we have proved in the following theorem:
\begin{theorem}
The optimal portfolio $\Pi^*(t)$ for the insider portfolio problem \eqref{eq17} is given by
\begin{equation}
\Pi^*(t) = \int_{\mathbb{R}}\pi^*(t,y) \delta_Y(y) dy=\pi^*(t,Y),
\end{equation}
where
\begin{equation}
\pi^*(t,y)=\frac{k(t,y)}{x(t,y)\gamma_0(t,y)}
\end{equation}
with $x(t,y), k(t,y)$ given as the solution of the BSDE \eqref{eq37'} and $c=p(0,y)$ given by \eqref{eq39'}.
\end{theorem}
\subsection{The logarithmic utility case (Brownian-Poisson process)}


We now extend the financial applications in the previous sections to the  case with both a Brownian motion component $B(t)$ and a Poisson process $N(t)$ with intensity $\lambda>0$, here we have  $\tilde{N}(t)=N(t) -\lambda t$. Thus we consider a financial market where the unit price $S_0(t)$ of the risk free asset is
\begin{equation}\label{eq6.78}
    S_0(t)=1, \quad t\in[0, T]
\end{equation}
and the unit price process $S(t)$ of the risky asset is given by

\begin{align}\label{eq6.79}
\begin{cases}
 dS(t) &= S(t) [b_0(t,Y) dt + \sigma_0(t,Y) dB(t) + \gamma_0(t,Y) d\tilde{N}(t)]
 ; \quad t\in[0, T]\\
    S(0)&>0.
\end{cases}
\end{align}
Then the wealth process $X(t)=X^{\Pi}(t)$ associated to a portfolio $u(t)=\Pi(t)$, interpreted as the fraction of the wealth invested in the risky asset at time $t$, is given by
\begin{align}\label{eq6.80}
  \begin{cases}
dX(t)&=\Pi(t)X(t)[b_0(t,Y)dt+\sigma_0(t,Y)dB(t)+\gamma_0(t,Y)d\tilde{N}(t)]
; \quad t\in[0, T]\\
X(0)&=x_0>0.
\end{cases}
  \end{align}

Let $U$ be a given utility function. We want to find $\Pi^{\ast}\in \mathcal{A}$ such that
\begin{equation}\label{eq6.81}
J(\Pi^*) =\sup_{\Pi\in\mathcal{A}} J(\Pi),
\end{equation}
where
\begin{equation}\label{eq6.82}
J(\Pi):=  \mathbb{E}[U(X^\Pi(T))].
\end{equation}
Note that, in terms of our process $x(t,y)$ we have
\begin{align}\label{eq6.83}
\begin{cases}
dx(t,y)=\pi(t,y)x(t,y)[b_0(t,y)dt+\sigma_0(t,y)dB(t)+ \gamma_0(t,y)d\tilde{N}(t)]
; \quad t\in[0, T]\\
x(0,y)=x_0(y)>0,
\end{cases}
\end{align}
which has the solution
\begin{align}\label{eq6.84}
x(t,y)& = x_0(y) \exp \big( \int_0^t \{\pi(s,y) b_0(s,y)-\frac{1}{2}\pi^2(s,y)\sigma_0^2(s,y) \nonumber\\
&+\lambda\ln(1+\pi(s,y)\gamma_0(s,y)) -\lambda \pi(s,y) \gamma_0(s,y)\} ds \nonumber\\
&+\int_0^t \pi(s,y) \sigma_0(s,y)dB(s) \nonumber\\
&+ \int_0^t  \ln(1+ \pi(s,y) \gamma_0(s,y)) d\tilde{N}(s) \big).
\end{align}
The performance functional gets the form
\begin{equation}\label{eq6.85}
J(\pi)=\mathbb{E}[ U(x(T,y)) \mathbb{E}[\delta_Y(y)|\mathcal{F}_T]],
\end{equation}
where
\begin{equation}\label{eq6.86}
\Pi(t)=\pi(t,Y).
\end{equation}
In this case the Hamiltonian (\ref{eq11}) gets the form
\begin{equation}\label{eq6.87}
    H(t,x,y,\pi,p,q,r)=\pi x[b_0(t,y)p+\sigma_0(t,y)q+\lambda\gamma_0(t,y)r(y,1)],
\end{equation}
while the BSDE for the adjoint processes becomes
\begin{equation}\label{eq6.88}
    \left\{
  \begin{array}{l}
dp(t,y)=- \pi(t,y)  [b_0(t,y)p(t,y) +\sigma_0(t,y)q(t,y)+\lambda\gamma_0(t,y) r(t,y,1) ]dt\\
+q(t,y)dB(t)+r(t,y,1)d\tilde{N}(t); \quad t\geq 0\\
p(T,y)=U'(x(T,y))\mathbb{E}[\delta_Y(y)|\mathcal{F}_T]
  \end{array}
    \right.
\end{equation}

If $U(x)$ is the logarithmic utility, i.e.
$$U(x)=\ln x; x > 0,$$
then
\begin{eqnarray}\label{eq6.89}
   \mathbb{E}[\ln(x(T,y))\mathbb{E}[\delta_Y(y)|\mathcal{F}_T]]&=& \mathbb{E}\big[\{\int_0^T (\pi(t,y) b_0(t,y)-\frac{1}{2}\pi^2(t,y)\sigma_0^2(t,y) \nonumber \\
   &+&\lambda \ln(1+\pi(t,y)\gamma_0(t,y)) - \lambda \pi(t,y) \gamma_0(t,y))dt \nonumber  \\
   &+& \int_0^T \pi(t,y) \sigma_0(t,y)dB(t) \nonumber \\
   &+&  \int_0^T \ln(1+ \pi(t,y) \gamma_0(t,y)) d\tilde{N}(t)\}\mathbb{E}[\delta_Y(y)|\mathcal{F}_T]\big]
\end{eqnarray}

We now use the duality formulas, Theorem $7.11$ and Theorem $8.5$. This enables us to write (\ref{eq6.89}) as the expectation of a $ds$-integral and we get:
 \begin{eqnarray}\label{eq6.90}
   \mathbb{E}[\ln(x(T,y))\mathbb{E}[\delta_Y(y)|\mathcal{F}_T]]&=& \mathbb{E}\big[\int_0^T \{\pi(t,y) b_0(t,y)-\frac{1}{2}\pi^2(t,y)\sigma_0^2(t)+\lambda\ln(1+\pi(t,y)\gamma_0(t,y)) \nonumber \\
   &-& \lambda\pi(t,y) \gamma_0(t,y)\}dt\mathbb{E}[\delta_Y(y)|\mathcal{F}_T] \nonumber \\
   &+& \int_0^T\mathbb{E}[D_t\mathbb{E}[\delta_Y(y)|\mathcal{F}_T]|\mathcal{F}_t] \pi(t,y) \sigma_0(t,y)dt \\
   &+&  \int_0^T\mathbb{E}[D_{t,1}\mathbb{E}[\delta_Y(y)|\mathcal{F}_T]|\mathcal{F}_t] \ln(1+ \pi(t,y) \gamma_0(t,y))\lambda dt\big]
\end{eqnarray}

Note that
\begin{equation}
D_t\mathbb{E}[\delta_Y(y)|\mathcal{F}_T]=\mathbb{E}[D_t\delta_Y(y)|\mathcal{F}_T]
\end{equation}
and
\begin{equation}
D_{t,1}\mathbb{E}[\delta_Y(y)|\mathcal{F}_T]=\mathbb{E}[D_{t,1}\delta_Y(y)|\mathcal{F}_T].
\end{equation}

Therefore, if we substitute this in \eqref{eq6.90} and take for each $t$ the conditional expectation with respect to $\mathcal{F}_t$ of the integrand, we get
\begin{eqnarray}\label{eq6.93}
   \mathbb{E}[\ln(x(T,y))\mathbb{E}[\delta_Y(y)|\mathcal{F}_T]]&=& \mathbb{E}\big[\int_0^T \{\pi(t,y) b_0(t,y)-\frac{1}{2}\pi^2(t,y)\sigma_0^2(t,y)+\lambda\ln(1+\pi(t,y)\gamma_0(t,y)) \nonumber\\
   &-&\lambda \pi(t,y) \gamma_0(t,y)\}dt\mathbb{E}[\delta_Y(y)|\mathcal{F}_T] \nonumber \\
   &+& \int_0^T\mathbb{E}[D_t\delta_Y(y)|\mathcal{F}_t] \pi(t,y) \sigma_0(t,y)dt \nonumber \\
   &+&  \int_0^T\mathbb{E}[D_{t,1}\delta_Y(y)|\mathcal{F}_t] \ln(1+ \pi(t,y) \gamma_0(t,y)) \lambda dt\big]
\end{eqnarray}

 We can maximize this by maximizing the integrand with respect to $\pi(t,y)$ for each $t$ and $y$.
 Doing this we obtain that the optimal portfolio $\pi(t,y)$ for Problem (\ref{eq6.81})  is given implicitly as the solution $\pi(t,y)$ of the first order condition

 \begin{align}\label{eq6.94}
& [b_0(t,y)-\pi(t,y)\sigma_0^2(t,y)-\frac{\lambda\pi(t,y) \gamma_0^2(t,y)}{1+\pi(t,y)\gamma_0(t,y)}]\mathbb{E}[\delta_Y(y)|\mathcal{F}_t]
+\sigma_0(t,y)\mathbb{E}[D_{t}\delta_Y(y)|\mathcal{F}_t] \nonumber\\
&+\frac{\lambda\gamma_0(t,y)}{1+\pi(t,y)\gamma_0(t,y)}\mathbb{E}[D_{t,1}\delta_Y(y)|\mathcal{F}_t]= 0.
\end{align}

If we define
\begin{equation}\label{eq6.95}
\Phi (t,y) := \frac{\mathbb{E}[D_t \delta_Y(y)|\mathcal{F}_t]}{\mathbb{E}[\delta_Y(y)|\mathcal{F}_t]}
\end{equation}
and
\begin{equation}\label{eq6.96}
\Psi(t,y):= \frac{\mathbb{E}[D_{t,1}\delta_Y(y)|\mathcal{F}_t]}{\mathbb{E}[\delta_Y(y)|\mathcal{F}_t]}
\end{equation}

then \eqref{eq6.94} can be written

\begin{align}\label{eq6.97}
& b_0(t,y)-\pi(t,y)\sigma_0^2(t,y)- \frac{\lambda \pi(t,y) \gamma_0^2(t,y)}{1+\pi(t,y)\gamma_0(t,y)}\nonumber\\
&+\sigma_0(t,y)\Phi(t,y)+\frac{\lambda \gamma_0(t,y)}{1+\pi(t,y)\gamma_0(t,y)}\Psi(t,y)= 0.
\end{align}

Thus we have proved the following theorem:

\begin{theorem}
The optimal portfolio with respect to logarithmic utility for an insider in the market (\ref{eq6.78})-\eqref{eq6.79}and with the inside information (\ref{H_t}) is given implicitly as the solution $\Pi (t)= \Pi^*(t)$ of the equation
 \begin{align}\label{eq6.98}
& b_0(t,Y)-\Pi(t)\sigma_0^2(t,Y)- \frac{\lambda\Pi(t) \gamma_0^2(t,Y)}{1+\Pi(t)\gamma_0(t,Y)}\nonumber\\
&+\sigma_0(t,Y)\Phi(t,Y)+\frac{\lambda\gamma_0(t,Y)}{1+\Pi(t)\gamma_0(t,Y)}\Psi(t,Y)= 0,
\end{align}
provided that a solution exists.
\end{theorem}



The equation \eqref{eq6.98} for the optimal portfolio $\Pi(t)$ holds for a general insider random variable $Y$. In the case when $Y$ is of the form \eqref{Brow-Poiss}, then we can substitute  (\ref{eq2.22}), (\ref{eq2.22a}) and \eqref{eq2.23} in (\ref{eq6.95}) and (\ref{eq6.96}), and get a more explicit equation as follows:\\
\begin{theorem}
Suppose $Y$ is as in \eqref{Brow-Poiss}. Then the processes $\Phi(t,y)$ and $\Psi(t,y)$ in the equation \eqref{eq6.98} for the optimal portfolio $\pi(t,y)$  have the following expressions:
\begin{align}
&\Phi(t,y)= \frac{i\beta \int_{\mathbb{R}} F(t,x,y) x dx}{\int_{\mathbb{R}} F(t,x,y) dx}\\
&\Psi(t,y) = \frac{\int_{\mathbb{R}} F(t,x,y) (e^{ix}-1)dx}{\int_{\mathbb{R}} F(t,x,y) dx}\\
\end{align}
where
\begin{align}
&F(t,x,y)=\frac{1}{2\pi}\exp\big[ix\tilde{N}(t) +ix\beta B(t)\nonumber\\
&   +\lambda(T_0-t)(e^{ix}-1-ix)-\frac{1}{2}x^2\beta^2(T_0-t)-ixy\big].
   \end{align}
   \end{theorem}


\subsection{The general It\^{o}-L\'{e}vy process case}
We now extend the financial applications in the previous sections to the general case with both a Brownian motion component $B(t)$ and a compensated Poisson random measure component $\tilde{N}(dt,d\zeta)=N(dt,d\zeta) -\nu(d\zeta)dt$, as in Section 2. Thus we consider a financial market where the unit price $S_0(t)$ of the risk free asset is
\begin{equation}\label{eq7.1}
    S_0(t)=1, \quad t\in[0, T]
\end{equation}
and the unit price process $S(t)$ of the risky asset is given by

\begin{align}\label{eq7.2}
\begin{cases}
 dS(t) &= S(t) [b_0(t,Y) dt + \sigma_0(t,Y) dB(t) +\int_{\mathbb{R}} \gamma_0(t,Y,\zeta) \tilde{N}(dt,d\zeta)]
 ; \quad t\in[0, T]\\
    S(0)&>0.
\end{cases}
\end{align}
Then the wealth process $X(t)=X^{\Pi}(t)$ associated to a portfolio $u(t)=\Pi(t)$, interpreted as the fraction of the wealth invested in the risky asset at time $t$, is given by
\begin{align}\label{eq7.3}
  \begin{cases}
dX(t)&=\Pi(t)X(t)[b_0(t,Y)dt+\sigma_0(t,Y)dB(t)+\int_{\mathbb{R}}\gamma_0(t,Y,\zeta)\tilde{N}(dt,d\zeta)]
; \quad t\in[0, T]\\
X(0)&=x_0>0.
\end{cases}
  \end{align}

Let $U$ be a given utility function. We want to find $\Pi^{\ast}\in \mathcal{A}$ such that
\begin{equation}\label{eq17'}
J(\Pi^*) =\sup_{\Pi\in\mathcal{A}} J(\Pi),
\end{equation}
where
\begin{equation}\label{eq18'}
J(\Pi):=  \mathbb{E}[U(X^\Pi(T))].
\end{equation}
Note that, in terms of our process $x(t,y)$ we have
\begin{align}\label{Wealth7}
\begin{cases}
dx(t,y)=\pi(t,y)x(t,y)[b_0(t,y)dt+\sigma_0(t,y)dB(t)+\int_\mathbb{R} \gamma_0(t,y,\zeta)\tilde{N}(dt,d\zeta)]
; \quad t\in[0, T]\\
x(0,y)=x_0(y)>0,
\end{cases}
\end{align}
which has the solution
\begin{align}\label{eq7.7}
x(t,y)& = x_0(y) \exp \big( \int_0^t \{\pi(s,y) b_0(s,y)-\frac{1}{2}\pi^2(s,y)\sigma_0^2(s,y) \nonumber\\
&+\int_0^t \int_{\mathbb{R}} [\ln(1+\pi(s,y)\gamma_0(s,y,\zeta)) - \pi(s,y) \gamma_0(s,y,\zeta)] \nu(d\zeta)\} ds \nonumber\\
&+\int_0^t \pi(s,y) \sigma_0(s,y)dB(s) \nonumber\\
&+ \int_0^t \int_{\mathbb{R}} \ln(1+ \pi(s,y) \gamma_0(s,y,\zeta)) \tilde{N}(ds,d\zeta) \big).
\end{align}
The performance functional gets the form
\begin{equation}\label{eq7.8}
J(\pi)=\mathbb{E}[ U(x(T,y)) \mathbb{E}[\delta_Y(y)|\mathcal{F}_T]],
\end{equation}
where
\begin{equation}\label{eq7.9}
\Pi(t)=\pi(t,Y).
\end{equation}
In this case the Hamiltonian (\ref{eq11}) gets the form
\begin{equation}\label{eq7.10}
    H(t,x,y,\pi,p,q,r)=\pi x[b_0(t,y)p+\sigma_0(t,y)q+\int_{\mathbb{R}}\gamma_0(t,y,\zeta)r(y,\zeta)\nu(d\zeta)],
\end{equation}
while the BSDE for the adjoint processes becomes
\begin{equation}\label{eq7.11}
    \left\{
  \begin{array}{l}
dp(t,y)=- \pi(t,y)  [b_0(t,y)p(t,y) +\sigma_0(t,y)q(t,y)+\int_{\mathbb{R}} \gamma_0(t,y,\zeta) r(t,y,\zeta)\nu(d\zeta) ]dt\\
+q(t,y)dB(t)+\int_{\mathbb{R}}r(t,y,\zeta)\tilde{N}(dt,d\zeta); \quad t\geq 0\\
p(T,y)=U'(x(T,y))\mathbb{E}[\delta_Y(y)|\mathcal{F}_T]
  \end{array}
    \right.
\end{equation}

If $U(x)$ is the logarithmic utility, i.e.
$$U(x)=\ln x; x > 0,$$
then
\begin{eqnarray}\label{eq7.12}
   \mathbb{E}[\ln(x(T,y))\mathbb{E}[\delta_Y(y)|\mathcal{F}_T]]&=& \mathbb{E}\big[\{\int_0^T \{\pi(t,y) b_0(t,y)-\frac{1}{2}\pi^2(t,y)\sigma_0^2(t,y) \nonumber \\
   &+&\int_0^T \int_{\mathbb{R}} [\ln(1+\pi(t,y)\gamma_0(t,y,\zeta)) - \pi(t,y) \gamma_0(t,y,\zeta)] \nu(d\zeta)\} dt \nonumber  \\
   &+& \int_0^T \pi(t,y) \sigma_0(t,y)dB(t) \nonumber \\
   &+&  \int_0^T \int_{\mathbb{R}} \ln(1+ \pi(t,y) \gamma_0(t,y,\zeta)) \tilde{N}(dt,d\zeta)\}\mathbb{E}[\delta_Y(y)|\mathcal{F}_T]\big]
\end{eqnarray}

We now use the duality formulas, Theorem $7.11$ and Theorem $8.5$. This enables us to write (\ref{eq7.12}) as the expectation of a $ds$-integral and we get:
 \begin{eqnarray}\label{eq7.13}
   \mathbb{E}[\ln(x(T,y))\mathbb{E}[\delta_Y(y)|\mathcal{F}_T]]&=& \mathbb{E}\big[\int_0^T \{\pi(t,y) b_0(t,y)-\frac{1}{2}\pi^2(t,y)\sigma_0^2(t,y)\}dt\mathbb{E}[\delta_Y(y)|\mathcal{F}_T] \nonumber \\
   &+&\int_0^T \int_{\mathbb{R}} [\ln(1+\pi(t,y)\gamma_0(t,y,\zeta)) - \pi(t,y) \gamma_0(t,y,\zeta)] \nu(d\zeta) dt\mathbb{E}[\delta_Y(y)|\mathcal{F}_T] \nonumber  \\
   &+& \int_0^T\mathbb{E}[D_t\mathbb{E}[\delta_Y(y)|\mathcal{F}_T]|\mathcal{F}_t] \pi(t,y) \sigma_0(t,y)dt \nonumber \\
   &+&  \int_0^T \int_{\mathbb{R}}\mathbb{E}[D_{t,\zeta}\mathbb{E}[\delta_Y(y)|\mathcal{F}_T]|\mathcal{F}_t] \ln(1+ \pi(t,y) \gamma_0(t,y,\zeta)) \nu(d\zeta)dt\big]
\end{eqnarray}

Note that
\begin{equation}
D_t\mathbb{E}[\delta_Y(y)|\mathcal{F}_T]=\mathbb{E}[D_t\delta_Y(y)|\mathcal{F}_T]
\end{equation}
and
\begin{equation}
D_{t,\zeta}\mathbb{E}[\delta_Y(y)|\mathcal{F}_T]=\mathbb{E}[D_{t,\zeta}\delta_Y(y)|\mathcal{F}_T].
\end{equation}

Therefore, if we substitute this in \eqref{eq7.13} and take for each $t$ the conditional expectation with respect to $\mathcal{F}_t$ of the integrand, we get
\begin{eqnarray}\label{eq7.13a}
   \mathbb{E}[\ln(x(T,y))\mathbb{E}[\delta_Y(y)|\mathcal{F}_T]]&=& \mathbb{E}\big[\int_0^T \{\pi(t,y) b_0(t,y)-\frac{1}{2}\pi^2(t,y)\sigma_0^2(t,y)\}dt\mathbb{E}[\delta_Y(y)|\mathcal{F}_T] \nonumber \\
   &+&\int_0^T \int_{\mathbb{R}} [\ln(1+\pi(t,y)\gamma_0(t,y,\zeta)) - \pi(t,y) \gamma_0(t,y,\zeta)] \nu(d\zeta) dt\mathbb{E}[\delta_Y(y)|\mathcal{F}_T] \nonumber  \\
   &+& \int_0^T\mathbb{E}[D_t\delta_Y(y)|\mathcal{F}_t] \pi(t,y) \sigma_0(t,y)dt \nonumber \\
   &+&  \int_0^T \int_{\mathbb{R}}\mathbb{E}[D_{t,\zeta}\delta_Y(y)|\mathcal{F}_t] \ln(1+ \pi(t,y) \gamma_0(t,y,\zeta)) \nu(d\zeta)dt\big]
\end{eqnarray}

 We can maximize this by maximizing the integrand with respect to $\pi(t,y)$ for each $t$ and $y$.
 Doing this we obtain that the optimal portfolio $\pi(t,y)$ for Problem (\ref{eq17'})  is given implicitly as the solution $\pi(t,y)$ of the first order condition

 \begin{align}\label{eq7.14}
& [b_0(t,y)-\pi(t,y)\sigma_0^2(t,y)]\mathbb{E}[\delta_Y(y)|\mathcal{F}_t]+\sigma_0(t,y)\mathbb{E}[D_{t}\delta_Y(y)|\mathcal{F}_t] \nonumber\\
& - \int_{\mathbb{R}} \frac{\pi(t,y) \gamma_0^2(t,y,\zeta)}{1+\pi(t,y)\gamma_0(t,y,\zeta)}[\mathbb{E}[\delta_Y(y)|\mathcal{F}_t] \nonumber\\
&+ \int_{\mathbb{R}}\frac{\gamma_0(t,y,\zeta)}{1+\pi(t,y)\gamma_0(t,y,\zeta)}\mathbb{E}[D_{t,\zeta}\delta_Y(y)|\mathcal{F}_t]= 0.
\end{align}

If we define
\begin{equation}\label{eq8.18}
\Phi (t,y) := \frac{\mathbb{E}[D_t \delta_Y(y)|\mathcal{F}_t]}{\mathbb{E}[\delta_Y(y)|\mathcal{F}_t]}
\end{equation}
and
\begin{equation}\label{eq8.19}
\Psi(t,\zeta,y):= \frac{\mathbb{E}[D_{t,\zeta}\delta_Y(y)|\mathcal{F}_t]}{\mathbb{E}[\delta_Y(y)|\mathcal{F}_t]}
\end{equation}

then \eqref{eq7.14} can be written

\begin{align}\label{eq7.20}
& b_0(t,y)-\pi(t,y)\sigma_0^2(t,y)- \int_{\mathbb{R}} \frac{\pi(t,y) \gamma_0^2(t,y,\zeta)}{1+\pi(t,y)\gamma_0(t,y,\zeta)} \nu(d\zeta)\nonumber\\
&+\sigma_0(t,y)\Phi(t,y)+ \int_{\mathbb{R}}\frac{\gamma_0(t,y,\zeta)}{1+\pi(t,y)\gamma_0(t,y,\zeta)}\Psi(t,y,\zeta)\nu(d\zeta)= 0.
\end{align}

Thus we have proved the following theorem:

\begin{theorem}
The optimal portfolio with respect to logarithmic utility for an insider in the market (\ref{eq7.1})-\eqref{eq7.2}and with the inside information (\ref{H_t}) is given implicitly as the solution $\Pi (t)= \Pi^*(t)$ of the equation
 \begin{align}\label{eq7.15}
& b_0(t,Y)-\Pi(t)\sigma_0^2(t,Y)- \int_{\mathbb{R}} \frac{\Pi(t) \gamma_0^2(t,Y,\zeta)}{1+\Pi(t)\gamma_0(t,Y,\zeta)} \nu(d\zeta)\nonumber\\
&+\sigma_0(t,y)\Phi(t,Y)+ \int_{\mathbb{R}}\frac{\gamma_0(t,Y,\zeta)}{1+\Pi(t)\gamma_0(t,Y,\zeta)}\Psi(t,Y,\zeta)\nu(d\zeta)= 0,
\end{align}
provided that a solution exists.
\end{theorem}



The equation \eqref{eq7.15} for the optimal portfolio $\Pi(t)$ holds for a general insider random variable $Y$. In the case when $Y$ is of the form \eqref{eq2.5}, then we can substitute  (\ref{eq2.9}), (\ref{2.15}) and \eqref{eq2.13} in (\ref{eq8.18}) and (\ref{eq8.19}), and get a more explicit equation as follows:\\
\begin{theorem}
Suppose $Y$ is as in \eqref{eq2.5}. Then the processes $\Phi(t,y)$ and $\Psi(t,y,z)$ in the equation \eqref{eq7.15} for the optimal portfolio $\pi(t,y)$  have the following expressions:
\begin{align}
&\Phi(t,y)= \frac{i\beta(t) \int_{\mathbb{R}} F(t,x,y) x dx}{\int_{\mathbb{R}} F(t,x,y) dx}\\
&\Psi(t,y,z) = \frac{\int_{\mathbb{R}} F(t,x,y) (e^{ix\psi(t,z)}-1)dx}{\int_{\mathbb{R}} F(t,x,y) dx}\\
\end{align}
where
\begin{align}
&F(t,x,y)=\int_{\mathbb{R}}\exp\big[\int_0^t\int_{\mathbb{R}}ix\psi(s,\zeta)\tilde{N}(ds,d\zeta) +\int_0^t ix\beta(s)dB(s)\nonumber\\
&   +\int_t^{T_0}\int_{\mathbb{R}}(e^{ix\psi(s,\zeta)}-1-ix\psi(s,\zeta))\nu(d\zeta)ds-\int_t^{T_0}\frac{1}{2}x^2\beta^2(s)ds-ixy\big]dx.
   \end{align}
   \end{theorem}



\section{Appendix}
For the convenience of the reader, we give in this Appendix a brief survey of the main concepts and results from the theory of Hida-Malliavin calculus and white noise analysis needed in the previous sections. For more details see e.g. \cite{BBS}, \cite {DOP}, \cite{DMOP1}, \cite{HOUZ}, \cite{OR2} and the references therein. \\

\subsection{The White Noise Probability Space and Hida-Malliavin Calculus for Brownian motion}
Let $\mathcal{S}(\mathbb{R})$ be the Schwartz space consisting of all real-valued rapidly decreasing
functions $f$ on $\mathbb{R},$ i.e.,
\begin{equation}\label{rapidely decreasing}
    \lim_{|x|\rightarrow \infty}|x^nf^{(k)}(x)|=0, \quad \forall n,k\geq 0
\end{equation}
For instance $\mathcal{C}^{\infty}$ functions with compact support $e^{-x^2}, e^{-x^4}, ...$ are all
functions in $\mathcal{S}(\mathbb{R})$. For any $n, k \geq 0,$ define a norm $\|.\|_{n, k}$ on $\mathcal{S}(\mathbb{R})$ by
\begin{equation}\label{norm on S}
    \|f\|_{n, k}=\sup_{x\in\mathbb{R}}|x^nf^{(k)}(x)|.
\end{equation}
Then $(\mathcal{S}(\mathbb{R}, \{\|.\|_{n, k}, n, k \geq 0\})$ is a topological space. In fact, it is a nuclear space.

Let $\mathcal{S}'(\mathbb{R})$ be the dual space of $\mathcal{S}(\mathbb{R})$, the space of tempered distributions.
Let $B$ denote the family of all Borel subsets of $\mathcal{S}(\mathbb{R})$ equipped with the weak topology.
\begin{theorem}(Minlos)
Let $E$ be a nuclear space with dual space $E^{\ast}$.
 A complex-valued function $\phi$ on $E$ is the characteristic functional of a probability measure $\nu$  on $E^{\ast}$ ,i.e.,
 \begin{equation}\label{Minlos}
    \phi(y)=\int_{E^{\ast}}e^{i\langle x, y \rangle} d\nu(x), \quad y\in E
 \end{equation}
  if and only if it satisfies the following conditions:
 \begin{enumerate}
 \item $\phi(0)=1$,
 \item $\phi$ is positive definite,
 \item $\phi$ is continuous.
\end{enumerate}
\end{theorem}

\begin{remark}
The measure $\nu$ is uniquely determined by $\phi$. Observe that $\phi(0)=\nu(E^{\ast}).$
Thus when condition $(1)$ is not assumed, then we can only conclude that $\nu$ is a finite measure.
\end{remark}

Let $\phi$ be a function on $\mathcal{S}(\mathbb{R})$ given by
\begin{equation*}
    \phi(\xi)=\exp -\frac{1}{2}|\xi|^2, \quad \xi\in \mathcal{S}(\mathbb{R})
\end{equation*}
where $|.|$ is the $\mathbf{L}^2(\mathbb{R})$ norm. Then it is easy to check that conditions $(1)$ and $(2)$ are satisfied.
To check condition $(3)$ note that
\begin{eqnarray*}
  |\xi|^2 &=& \int_{\mathbb{R}}|\xi(x)|^2 dx \\
   &=& \int_{|x| < 1}|\xi(x)|^2 dx + \int_{|x| \geq 1}|\xi(x)|^2 dx  \\
   &\leq& 2\sup_{|x| < 1}|\xi(x)|^2 + \int_{|x| \geq 1}\frac{1}{x^2}|x\xi(x)|^2 dx \\
   &\leq& 2\|\xi\|^2_{0, 0} + \sup_{|y| \geq 1}|y\xi(y)|^2\int_{|x| \geq 1}\frac{1}{x^2} dx   \\
   &\leq& 2\|\xi\|^2_{0, 0} +2\|\xi\|^2_{1, 0}
\end{eqnarray*}
This shows that $\phi$ is continuous.
Therefore, by the Minlos theorem there
exists a unique probability measure $\mathbf{P}$ on $\mathcal{S}'(\mathbb{R})$ such that
\begin{definition}
The measure $\mathbf{P}$ is called the standard Gaussian measure on $\mathcal{S}'(\mathbb{R})$.
The probability space $(\mathcal{S}'(\mathbb{R}), \mathcal{B}, \mathbf{P})$ is called a white noise space.
\end{definition}
The Schwartz distribution theory on the space $\mathbb{R}$ concerns with the following Gel'fand triple
\begin{equation}\label{Gelfand1}
    \mathcal{S}(\mathbb{R}) \subset \mathbf{L}^2(\mathbb{R}) \subset \mathcal{S}'(\mathbb{R})
\end{equation}

\subsection{The Wiener It\^{o} chaos expansion}
let the Hermite polynomials $h_n(x)$ be defined by
\begin{equation*}
    h_n(x)=(-1)^n e^{\frac{1}{2}x^2}\frac{d^n}{dx^n}(e^{-\frac{1}{2}x^2}) \quad n= 0, 1, 2, ...
\end{equation*}
Let $e_k$ be the kth Hermite function defined by
\begin{equation}\label{hermite function}
    e_k(x) := \pi^{-\frac{1}{4}}((k-1)!)^{-\frac{1}{2}}e^{-\frac{1}{2}x^2}h_{k-1}(\sqrt{2}x), \quad k = 1, 2, ...
\end{equation}
Then $\{e_k\}_{k\geq 1}$ constitutes an orthonormal basis for $\mathbf{L}^2(\mathbf{R})$ and $e_k \in \mathcal{S}(\mathbb{R})$ for
all $k$.
Define
\begin{equation}\label{theta_k}
    \theta_k(\omega) := \langle \omega, e_k\rangle = \int_{\mathbb{R}}e_k(x)dB(x,\omega), \quad \omega \in\Omega
\end{equation}
Let $\mathcal{J}$ denote the set of all finite multi-indices $\alpha = (\alpha_1, \alpha_2, . . . , \alpha_m), m =
1, 2, . . .,$ of non-negative integers $\alpha_i.$ If $\alpha = (\alpha_1, . . . , \alpha_m) \in \mathcal{J} , \alpha \neq 0,$ we put
\begin{equation}\label{H_alpha}
    H_{\alpha}(\omega) := \prod_{j=1}^mh_{\alpha_j}(\theta_j(\omega)), \quad \omega\in\Omega
\end{equation}
By a result of It\^{o} we have that
\begin{equation}\label{I_m}
    I_m(e^{\widehat{\otimes}\alpha}) = \prod_{j=1}h_{\alpha_j}(\theta_j)=H_{\alpha}.
\end{equation}
We set $H_0 := 1$. Here and in the sequel the functions $e_1, e_2, . . .$ are
defined in (\ref{hermite function}) and $\otimes$ and $\widehat{\otimes}$ denote the tensor product and the symmetrized
tensor product, respectively.\\
The family $\{H_{\alpha}\}_{\alpha\in \mathcal{J}}$ is an orthogonal basis for the Hilbert space $\mathbf{L}^2(\mathbf{P})$. In
fact, we have the following result.
\begin{theorem}\label{wiener ito1}The Wiener It\^{o} chaos expansion theorem. The family
$\{H_{\alpha}\}\alpha \in \mathcal{J}$ constitutes an orthogonal basis of $\mathbf{L}^2(\mathbf{P})$. More precisely, for all
$F_T$ -measurable $X \in\mathbf{L}^2(\mathbf{P})$  there exist (uniquely determined) numbers $c_{\alpha}\in\mathbb{R}$
such that
\begin{equation}\label{decomposition1}
    X =\sum_{\alpha\in\mathcal{J}}c_{\alpha}H_{\alpha} \in \mathbf{L}^2(\mathbf{P}).
\end{equation}
Moreover, we have
\begin{equation}\label{Norm1}
 \|X\|_{\mathbf{L}^2(\mathbf{P})}^2 = \sum_{\alpha\in\mathcal{J}}\alpha! c_{\alpha}^2.
\end{equation}
\end{theorem}
Let us compare Theorem (\ref{wiener ito1}) with the equivalent formulation of this theorem
in terms of iterated It\^{o} integrals. In fact, if $\psi(t_1, t_2, . . . , t_n)$
is a real symmetric function in its $n$ variables $t_1, . . . , t_n$ and $\psi\in \mathbf{L}^2(\mathbb{R}^n),$
that is,
\begin{equation*}
    \|\psi\|_{\mathbf{L}^2(\mathbb{R}^n)} := \big[\int_{\mathbb{R}^n}|\psi(t_1, t_2, ..., t_n)|^2 dt_1dt_2...dt_n\big]^\frac{1}{2}<\infty
\end{equation*}
then its $n$-tuple It\^{o} integral is defined by
\begin{eqnarray*}
I_n(\psi) &:=& \int_{\mathbb{R}^n}\psi dB^{\otimes n} \\
   &=& n! \int_{-\infty}^{\infty}\int_{-\infty}^{t_n}\int_{-\infty}^{t_{n-1}} ...\int_{-\infty}^{t_2}\psi(t_1, t_2, ..., t_n)dB(t_1)dB(t_2)...dB(t_n),
\end{eqnarray*}
where the integral on the right-hand side consists of $n$ iterated It\^{o} integrals.
Note that the integrand at each step is adapted to the filtration $\mathbb{F}$. Applying
the It\^{o} isometry $n$ times we see that
\begin{equation*}
    \mathbb{E}\big[(\int_{\mathbb{R}^n}\psi dB^{\otimes n})^2\big]=n!\|\psi\|^2_{\mathbf{L}^2(\mathbb{R}^n)}.
\end{equation*}
For $n = 0$ we adopt the convention that
\begin{equation*}
    I_0(\psi) := \int_{\mathbb{R}^0}\psi dB^{\otimes 0}=\psi=\|\psi\|_{\mathbf{L}^2(\mathbb{R}^0)},
\end{equation*}
for $\psi$ constant. Let $\widetilde{L}^2(\mathbb{R}^n)$ denote the set of symmetric real functions on $\mathbb{R}^n$,
which are square integrable with respect to Lebesque measure.
Then we have the following result
\begin{theorem}\label{wiener ito2} The Wiener It\^{o} chaos expansion theorem. For all $\mathcal{F}_t$-
measurable $X\in\mathbf{L}^2(\mathbf{P})$ there exist (uniquely determined) functions $f_n\in\widetilde{\mathbf{L}}^2(\mathbb{R}^n)$ such that
\begin{equation}\label{decomposition2}
    X=\sum_{n=0}^{\infty}\int_{\mathbb{R}^n}f_ndB^{\otimes n}=\sum_{n=0}^{\infty}I_n(f_n) \in \mathbf{L}^2(\mathbf{P})
\end{equation}
Moreover, we have the isometry
\begin{equation}\label{isometry2}
    \|X\|^2_{\mathbf{L}^2(\mathbf{P}) }= \sum_{n=0}^{\infty}n!\|f_n\|^2_{\mathbf{L}^2(\mathbb{R}^n)}
\end{equation}
\end{theorem}
The connection between these two expansions in Theorem (\ref{wiener ito1})
and Theorem (\ref{wiener ito2}) is given by
\begin{equation*}
    f_n=\sum_{\alpha\in\mathcal{J},|\alpha|=n}c_{\alpha}e_1^{\otimes \alpha_1}\widehat{\otimes}e_2^{\otimes \alpha_2}\widehat{\otimes}...\widehat{\otimes}e_m^{\otimes \alpha_m}, \quad n=0, 1, 2, ...
\end{equation*}
where $|\alpha| = \alpha_1+ \alpha_2...+ \alpha_m$ for $\alpha = (\alpha_1, ..., \alpha_m) \in\mathcal{J}, m = 1, 2, ...$ Recall that
the functions $e_1, e_2, ...$ are defined in (\ref{hermite function}) and $\otimes$ and $\widehat{\otimes}$ denote the tensor
product and the symmetrized tensor product, respectively.
Note that since $H_{\alpha} = I_m(e^{\widehat{\otimes}\alpha}),$ for $\alpha\in\mathcal{J}, |\alpha| = m,$ we get that
\begin{equation}\label{identit}
    m!\|e^{\widehat{\otimes}\alpha}\|^2_{\mathbf{L}^2(\mathbb{R}^m)} = \alpha!
\end{equation}
by combining (\ref{Norm1}) and (\ref{isometry2}) for $X = X_\alpha$.

Analogous to the test functions $\mathcal{S}(\mathbb{R})$  and the tempered distributions $\mathcal{S}'(\mathbb{R})$
on the real line $\mathbb{R},$ there is a useful space of stochastic test functions $(\mathcal{S})$ and
a space of stochastic distributions $(\mathcal{S}')$ on the white noise probability space.\\
In the following we use the notation
\begin{equation}\label{notation}
    (2\mathbb{N})^\alpha=\prod_{j=1}^{m}(2j)^{\alpha j}
\end{equation}
\subsection{The Kondratiev Spaces $(\mathcal{S})_1, (\mathcal{S})_{-1}$
and the Hida Spaces $(\mathcal{S})$ and $(\mathcal{S})^{\ast}$ }

\begin{definition}
Let $\rho$ be a constant in $[0,1]$.
\begin{enumerate}
\item Let $k\in\mathbb{R}$. We say that $f = \sum_{\alpha\in\mathcal{J}}a_{\alpha}H_{\alpha} \in \mathbf{L}^2(\mathbf{P})$ belongs to the Kondratiev test function Hilbert space $(\mathcal{S})_{k,\rho}$ if
    \begin{equation}\label{(S)_k norm}
        \|f\|^2_{k,\rho}:=\sum_{\alpha\in\mathcal{J}}a_{\alpha}^2 (\alpha !)^{1+\rho}(2\mathbb{N})^{\alpha k}<\infty.
    \end{equation}
We define the Kondratiev test function space $(\mathcal{S})_\rho$ as the space
\begin{equation*}
    (\mathcal{S})_\rho = \bigcap_{k\in\mathbb{R}}(\mathcal{S})_{k,\rho}
\end{equation*}
equipped with the projective topology, that is, $f_n\rightarrow f, n\rightarrow \infty,$ in $(\mathcal{S})_\rho$ if and only if $\|f_n-f\|_{k,\rho}\rightarrow 0,$ $n\rightarrow\infty,$ for all $k$.
\item Let $q\in\mathbb{R}$. We say that the formal sum $F = \sum_{\alpha\in\mathcal{J}}b_{\alpha}H_{\alpha}$ belongs to the Kondratiev stochastic distribution space $(\mathcal{S})_{-q,-\rho}$ if
     \begin{equation}\label{(S)_{-q} norm}
        \|f\|^2_{-q,-\rho}:=\sum_{\alpha\in\mathcal{J}}b_{\alpha}^2(\alpha!)^{1-\rho}(2\mathbb{N})^{-\alpha q}<\infty.
    \end{equation}
We define the Kondratiev distribution space $(\mathcal{S})_{-\rho}$ by
\begin{equation*}
    (\mathcal{S})_{-\rho} = \bigcup_{q\in\mathbb{R}}(\mathcal{S})_{-q,-\rho}
\end{equation*}
equipped with the inductive topology, that is, $F_n\rightarrow F, n\rightarrow \infty,$ in $(\mathcal{S})_{-\rho}$ if and only if there exists $q$ such that $\|F_n-F\|_{-q,-\rho}\rightarrow 0, n\rightarrow \infty.$
\item
If $\rho=0$ we write
\begin{equation}
(\mathcal{S})_0 = (\mathcal{S}) \text{ and } (\mathcal{S})_{-0} = (\mathcal{S})^{\ast}.
\end{equation}
These spaces are called the \emph{Hida test function space} and \emph{the Hida distribution space}, respectively.
\item If $F = \sum_{\alpha\in\mathcal{J}}b_{\alpha}H_{\alpha}$ in $(\mathcal{S})_{-1}$, we define the  generalized expectation $\mathbb{E}[F]$ of $F$ by
\begin{equation}\label{expectation}
    \mathbb{E}[F]=b_0.
\end{equation}
(Note that if $F \in\mathbf{L}^2(\mathbf{P})$, then the generalized expectation coincides with the usual expectation, since $\mathbb{E}[H_{\alpha}]=0$ for all $\alpha \neq 0$).
\end{enumerate}
\end{definition}
Note that $(\mathcal{S})_{-1}$ is the dual of $(\mathcal{S})_1$ and $(\mathcal{S})^{\ast}$ is the dual of $(\mathcal{S})$. The action of $F = \sum_{\alpha\in\mathcal{J}}b_{\alpha}H_{\alpha} \in (\mathcal{S})_{-1}$ on $f = \sum_{\alpha\in\mathcal{J}}a_{\alpha}H_{\alpha} \in (\mathcal{S})_1$ is given by
\begin{equation*}
    \langle F, f\rangle = \sum_{\alpha}\alpha! a_{\alpha}b_{\alpha}.
\end{equation*}
We have the inclusion
\begin{equation*}
   (\mathcal{S})_1\subset (\mathcal{S})\subset \mathbf{L}^2(\mathbf{P})\subset (\mathcal{S})^{\ast}\subset(\mathcal{S})_{-1}.
\end{equation*}
\subsection{The Spaces $\mathcal{G}$ and $\mathcal{G}^{\ast}$.}
We now introduce another pair of dual spaces, $\mathcal{G}$ and $\mathcal{G}^{\ast}$, which is sometimes
useful.
\begin{definition}
\begin{enumerate}
\item Let $\lambda \in \mathbb{R}$. Then the space $\mathcal{G}_{\lambda}$ consists of all formal
expansions
\begin{equation}\label{XinG}
    X=\sum_{n=0}^{\infty}\int_{\mathbb{R}^n}f_ndB^{\otimes n}
\end{equation}
such that
\begin{equation}\label{norm G_lambda}
    \|X\|_{\mathcal{G}_{\lambda}}=(\sum_{n=0}^{\infty}n! e^{2\lambda n}\|f_n\|^2_{\mathbf{L}^2(\mathbb{R}^n)})^{\frac{1}{2}}
\end{equation}
For each $\lambda \in \mathbb{R}$, the space $\mathcal{G}_{\lambda}$ is a Hilbert space with inner product
\begin{equation}\label{inner product}
    (X, Y)_{\mathcal{G}_{\lambda}}=\sum_{n=0}^{\infty}n! e^{2\lambda n}(f_n, g_n)_{\mathbf{L}^2(\mathbb{R}^n)}
\end{equation}
for every
\begin{equation*}
    X = \sum_{n=0}^{\infty}\int_{\mathbb{R}^n}f_ndB^{\otimes n}, \quad Y = \sum_{n=0}^{\infty}\int_{\mathbb{R}^n}g_ndB^{\otimes n}.
\end{equation*}
Note that $\lambda_1 \leq \lambda_2$ implies $\mathcal{G}_2 \subset \mathcal{G}_1$ . Define
\begin{equation}\label{G}
    \mathcal{G} = \bigcap_{\lambda \in\mathbb{R}} \mathcal{G}_{\lambda} = \bigcap_{\lambda > 0} \mathcal{G}_{\lambda},
\end{equation}
with projective limit topology.
\item $\mathcal{G}^{\ast}$ is defined to be the dual of $\mathcal{G}$. Hence
\begin{equation}\label{G^ast}
    \mathcal{G}^{\ast} = \bigcup_{\lambda \in\mathbb{R}} \mathcal{G}_{\lambda} = \bigcup _{\lambda < 0} \mathcal{G}_{\lambda},
\end{equation}
with inductive limit topology.
\end{enumerate}
\end{definition}

\begin{remark}
Note that an element $Y \in\mathcal{G}^{\ast}$ can be represented as a formal sum
\begin{equation}\label{YinG_ast}
    Y = \sum_{n=0}^{\infty}\int_{\mathbb{R}^n}g_ndB^{\otimes n}
\end{equation}
where $g_n\in \widetilde{\mathcal{L}}^ 2(\mathbb{R}^n)$ and $\|Y\|_{\mathcal{G}_{\lambda}}<\infty$ for some $\lambda \in \mathbb{R}$, while an $X \in \mathcal{G}$ satisfies $\|Y\|_{\mathcal{G}_{\lambda}}<\infty$ for all $\lambda \in \mathbb{R}$.
\end{remark}
If $X \in\mathcal{G}$ and $Y\in \mathcal{G}^{\ast}$. have the representations (\ref{XinG}) and (\ref{YinG_ast}), respectively,
then the action of $Y$ on $X$, $\langle Y, X \rangle$, is given by
\begin{equation}\label{duality G^ast,G }
    \langle Y, X \rangle = \sum_{n=0}^{\infty}n! (f_n, g_n)_{\mathcal{L}^ 2(\mathbb{R}^n)}.
\end{equation}
One can show that
\begin{equation*}
    (\mathcal{S})\subset\mathcal{G}\subset\mathbf{L}^2(\mathbf{P})\subset\mathcal{G}^{\ast}\subset(\mathcal{S})^{\ast}.
\end{equation*}
Finally, we note that, since
\begin{equation}\label{H_alpha}
    H_{\alpha} = \int_{\mathbb{R}^n}e^{\widehat{\otimes}\alpha}dB^{\otimes n} = I_n(e^{\otimes \alpha}),
\end{equation}
with $\alpha \in\mathcal{J}, |\alpha|=n,$ we get
\begin{equation*}
    \|H_{\alpha}\|^2_{\mathcal{G}_{\lambda}} = n! e^{2\lambda n} \|e^{\otimes \alpha}\|^2_{\mathcal{L}^ 2(\mathbb{R}^n)} = \alpha! e^{2\lambda n},
\end{equation*}
by (\ref{identit}). Therefore, for $F = \sum_{\alpha\in\mathcal{J}}c_{\alpha}H_{\alpha} \in \mathcal{G}^{\ast}$, we have
\begin{equation*}
    \|F_{\alpha}\|^2_{\mathcal{G}_{\lambda}} = \sum_{\alpha\in\mathcal{J}}c_{\alpha}^2 \alpha! e^{2\lambda |\alpha|}
\end{equation*}
for some $\lambda \in \mathbb{R}$.

\subsection{The Hida-Malliavin derivative}
\begin{definition}
\begin{enumerate}
\item Let $F \in\mathbf{L}^2(\mathbf{P})$ and let $h \in\mathbf{L}^2(\mathbb{R})$ be deterministic. Then the directional derivative of $F$ in $(\mathcal{S})^{\ast}$ in the direction $h$ is defined by
    \begin{equation}\label{def1directional deriv}
 D_hF (X):=\lim_{\epsilon \rightarrow 0}\frac{1}{\epsilon}[F(X+\epsilon h)-F(X)].
\end{equation}
whenever the limit exists in $(\mathcal{S})^{\ast}$.
\item Suppose there exists a function $\psi : \mathbb{R}\rightarrow (\mathcal{S})^{\ast}$ such that
$\int_{\mathbb{R}} \psi(t)h(t)dt$ converge in $(\mathcal{S})^{\ast}$ and
\begin{equation}\label{def2directional deriv}
    D_hF=\int_{\mathbb{R}} \psi(t)h(t)dt, \quad \text{for all } h \in \mathbf{L}^2(\mathbb{R}),
\end{equation}

then we say that $F$ is Hida-Malliavin differentiable in $(\mathcal{S})^{\ast}$ and we write
\begin{equation*}
    \psi(t)=D_tF, \quad t\in \mathbb{R}.
\end{equation*}
We call $D_tF$ the Hida-Malliavin derivative in $(\mathcal{S})^{\ast}$ or the stochastic gradient of $F$ at $t$.
\item More generally, if $F = \sum_{\alpha\in\mathcal{J}}b_{\alpha}H_{\alpha} \in (\mathcal{S})_{-1}$ we define the Hida-Malliavin derivative of $F$ at $t$ by the expansion
\begin{equation}
D_t F=\sum_{\alpha\in\mathcal{J}}\sum_{k=1}^{\infty}b_{\alpha}\alpha_k e_k(t)H_{\alpha-\epsilon^{(k)}}
\end{equation}
\end{enumerate}
\end{definition}
In \cite{AaOPU} $D_t$ it was shown that this is an extension from the space $D_{1,2}$ to $(\mathcal{S})_{-1}$
where $D_{1,2}$ denotes the classical space of Hida-Malliavin differentiable $F_T$-
measurable random variables. The extension is such that for all $F \in \mathbf{L}^2(F_T , \mathbf{P})$,
the following holds:
\begin{theorem}
\begin{enumerate}
\item Let $F \in (\mathcal{S})_{-1}$. Then $D_tF\in (\mathcal{S})_{-1} \text{ for all } t.$
\item Let $F \in\mathcal{G}^{\ast}$. Then $D_tF\in\mathcal{G}^{\ast}$ and $\mathbb{E}[D_tF|\mathcal{F}_t]\in\mathcal{G}^{\ast}$.
\item Let $F \in \mathbf{L}^2(F_T , \mathbf{P})$. Then $D_tF \in (\mathcal{S})^{\ast}.$
\item The map $(t, \omega) \mapsto \mathbb{E}[D_tF|\mathcal{F}_t]$ belongs to $\mathbf{L}^2(F_T, \lambda \times \mathbf{P})$, where $\lambda$ denotes
the Lebesgue measure on $[0, T]$.
\item
Moreover, the following generalized Clark-Ocone theorem holds:
\begin{equation}\label{Clark-Ocone}
    F = \mathbb{E}[F] + \int_0^T \mathbb{E}[D_tF|\mathcal{F}_t] dB(t)
\end{equation}
Notice that by combining It\^{o}'s isometry with the Clark-Ocone theorem, we
obtain
\begin{equation}\label{isometryD_t}
    \mathbb{E}[\int_0^T \mathbb{E}[D_tF|\mathcal{F}_t]^2 dt ] = \mathbb{E}[(\int_0^T\mathbb{E}[D_tF|\mathcal{F}_t] dB(t))^2] = \mathbb{E}[(F^2-\mathbb{E}[F]^2)]
\end{equation}
\end{enumerate}
\end{theorem}
As observed in \cite{AO} and \cite{DMOR}, we can apply the generalized Clark-Ocone theorem
to show that:
\begin{theorem}
(Generalized duality formula)\\
 Let $F \in \mathbf{L}^2(\mathcal{F}_T, \mathbf{P})$ and let
$\phi(t) \in \mathbf{L}^2(\lambda\times\mathbf{P})$ be adapted. Then
\begin{equation}\label{duality D_t}
    \mathbb{E}[F \int_0^T \phi(t) dB(t)] = \mathbb{E}[\int_0^T\mathbb{E}[D_tF|\mathcal{F}_t] \phi(t) dt]
\end{equation}
\end{theorem}
\begin{theorem}{Ordinary chain rule.}
\begin{enumerate}
\item Let $F\in\mathbf{L}^2(\mathbf{P})$ such that $F=\int_{\mathbb{R}} f(t)dB(t)$, $f\in\mathbf{L}^2(\mathbb{R})$. Then $F$ is Hida-Malliavin differentiable and
\begin{equation}
    D_t \int_{\mathbb{R}}f(s)dB(s)= f(t), t- a.a.
\end{equation}
\item Let $F\in\mathbf{L}^2(\mathbf{P})$ be Hida-Malliavin differentiable in $\mathbf{L}^2(\mathbf{P})$ for a.a. $t$. Suppose
that $\phi\in \mathbf{C}^1(\mathbb{R})$ and $\phi'(F)D_tF\in\mathbf{L}^2(\mathbf{P}\times\lambda)$.
Then $\phi(F)$ is Hida-Malliavin differentiable and we have
\begin{equation}
     D_t(\phi(F))=\phi'(F)D_tF,
\end{equation}
\end{enumerate}
\end{theorem}
\subsection{The Wick Product}
In addition to a canonical vector space structure, the spaces $(\mathcal{S})$ and $(\mathcal{S})^{\ast}$. also
have a natural multiplication given by the Wick product.
\begin{definition}
Let $X=\sum_{\alpha\in\mathcal{J}}a_{\alpha}H_{\alpha}$ and $Y =\sum_{\beta\in\mathcal{J}}b_{\beta}H_{\beta}$ be two elements of $(\mathcal{S})^{\ast}$. Then we define the Wick product of $X$ and $Y$ by
\begin{equation*}
    X\diamond Y=\sum_{\alpha,\beta\in\mathcal{J}}a_{\alpha}b_{\beta}H_{\alpha+\beta}=\sum_{\gamma\in\mathcal{J}}(\sum_{\alpha+\beta=\gamma}a_{\alpha}b_{\beta})H_{\gamma}
\end{equation*}
\end{definition}
For example we have
\begin{equation}\label{}
    B(t)\diamond B(t)=B^2(t)-t
\end{equation}
and more generally
\begin{equation}\label{}
    (\int_{\mathbb{R}}\phi(s)dB(s))\diamond(\int_{\mathbb{R}}\psi(s)dB(s))=(\int_{\mathbb{R}}\phi(s)dB(s)).(\int_{\mathbb{R}}\psi(s)dB(s))-\int_{\mathbb{R}}\phi(s)\psi(s)ds
\end{equation}
for all $\phi, \psi\in\mathbf{L}^2(\mathbb{R})$.

We list some properties of the Wick product:
\begin{enumerate}
\item $X, Y\in(\mathcal{S})_1 \Rightarrow X \diamond Y \in(\mathcal{S})_1$.
\item $X, Y\in(\mathcal{S})_{-1} \Rightarrow X \diamond Y \in(\mathcal{S})_{-1}$.
\item $X, Y \in(\mathcal{S}) \Rightarrow X \diamond Y \in(\mathcal{S})$.
\item $X \diamond Y = Y \diamond X$.
\item $X \diamond (Y \diamond Z) = (X \diamond Y) \diamond Z$.
\item $X \diamond (Y + Z) = X \diamond Y + X \diamond Z$.
\item $I_n(f_n) \diamond I_m(g_m) = I_{n+m}(f_n\hat{\otimes}g_m)$
\end{enumerate}
In view of the properties $(1)$ and $(4)$ we can define the Wick powers $X^{\diamond n}$
$(n = 1, 2, ...)$ of $X \in(\mathcal{S})_{-1}$ as
\begin{equation*}
  X^{\diamond n} := X\diamond X \diamond ...\diamond X \text{ (n times) }.
  \end{equation*}

We put $X^{\diamond 0} := 1$. Similarly, we define the Wick exponential $\exp^{\diamond}X$ of $X\in(\mathcal{S})_{-1}$  by
\begin{equation}\label{}
  \exp^{\diamond}X :=\sum_{n=0}^{\infty}\frac{1}{n!}X^{\diamond n} \in (\mathcal{S})_{-1}
  \end{equation}
   In view of the aforementioned properties, we
have that
\begin{equation}\label{}
  (X + Y )^{\diamond2} = X^{\diamond2} + 2X \diamond Y + Y^{\diamond2}
  \end{equation}

and also
\begin{equation}\label{}
  \exp^{\diamond}(X + Y ) = \exp^{\diamond}X \diamond \exp^{\diamond}Y,
  \end{equation}
for $X, Y\in\mathcal{S}_{-1}.$
Let $\mathbb{E}[.]$ denote the generalized expectation. Then
we see that
\begin{equation}\label{}
    \mathbb{E}[X \diamond Y ] = \mathbb{E}[X]\mathbb{E}[Y ],
\end{equation}
for $X, Y \in(\mathcal{S})_{-1}$. Note that independence is not required for this identity to
hold.
By induction, it follows that
\begin{equation}\label{}
    \mathbb{E}[\exp^{\diamond}X] = \exp{\mathbb{E}[X]},
\end{equation}
for $X \in(\mathcal{S})_{-1}$.
\begin{theorem} \textbf{Wick chain rule.}
\begin{enumerate}
\item Let $F,G\in (\mathcal{S})_{-1}$. Then $F\diamond G\in (\mathcal{S})_{-1}$ and
\begin{equation}
     D_t(F \diamond G) = F \diamond D_tG + D_tF \diamond G, t\in\mathbb{R}.
\end{equation}
\item Let $F \in(\mathcal{S})_{-1}$. Then
\begin{equation}
     D_t(F^{\diamond n}) = nF^{\diamond(n-1)} \diamond D_tF \quad (n = 1, 2, ...).
\end{equation}
\item
Let $F\in (\mathcal{S})_{-1}$. Then
\begin{equation}
     \exp^{\diamond}F =\sum_{n=0}^{\infty}\frac{1}{n!}F^{\diamond n} \in (\mathcal{S})_{-1}
\end{equation}

and
\begin{equation}
     D_t \exp^{\diamond}F = \exp^{\diamond}F \diamond D_tF.
\end{equation}
\end{enumerate}
\end{theorem}
\subsection {Conditional expectation}

If $F= \sum_{\alpha\in\mathcal{J}}b_{\alpha}H_{\alpha} \in (\mathcal{S})_{-1},$ we define its conditional expectation by the expansion
 \begin{equation}
\mathbb{E}[F|\mathcal{F}_t] :=  \sum_{\alpha\in\mathcal{J}}b_{\alpha}\mathbb{E}[H_{\alpha}| \mathcal{F}_t]
\end{equation}
Then the following holds:
\begin{enumerate}
\item
If $F,G \in (\mathcal{S})_{-1},$ then $\mathbb{E}[(F \diamond G) | \mathcal{F}_t] \in (\mathcal{S})_{-1}$ and \\
$\mathbb{E}[(F \diamond G) | \mathcal{F}_t] = \mathbb{E}[F|\mathcal{F}_t] \diamond \mathbb{E}[G|\mathcal{F}_t]$
\item
If $F,G \in (\mathcal{G})^{\ast},$ then $\mathbb{E}[(F \diamond G) | \mathcal{F}_t] \in (\mathcal{G})^{\ast}.$
\end{enumerate}
\subsection{The forward integral with respect to Brownian motion}
The forward integral with respect to Brownian motion was first defined in
the seminal paper \cite{RV} and further studied in \cite{RV1}, \cite{RV2}. This integral was
introduced in the modeling of insider trading in \cite{BO} and then applied by several authors
in questions related to insider trading and stochastic control with
advanced information (see, e.g., \cite{DMOP2}).
\begin{definition}
We say that a stochastic process $\phi = \phi(t), t\in[0, T ]$, is
forward integrable (in the weak sense) over the interval $[0, T ]$ with respect to
$B$ if there exists a process $I = I(t), t\in[0, T ]$, such that
\begin{equation}
\sup_{t\in[0,T ]}|\int_0^t\phi(s)\frac{B(s+\epsilon)-B(s)}{\epsilon}ds-I(t)|\rightarrow 0, \quad \epsilon\rightarrow0^+
\end{equation}
in probability. In this case we write
\begin{equation}
I(t) :=\int_0^t\phi(s)d^-B(s), t\in[0, T ],
\end{equation}
and call $I(t)$ the forward integral of $\phi$ with respect to $B$ on $[0, t]$.
\end{definition}
The following results give a more intuitive interpretation of the forward integral
as a limit of Riemann sums.
\begin{lemma}
Suppose $\phi$ is c\`{a}gl\`{a}d and forward integrable. Then
\begin{equation}
\int_0^T\phi(s)d^-B(s) = \lim _{\triangle t\rightarrow0}\sum_{j=1}^{J_n}\phi(t_{j-1})(B(t_j)-B(t_{j-1}))
\end{equation}
with convergence in probability. Here the limit is taken over the partitions \\
$0 =t_0 < t_1 < ... < t_{J_n}= T$ of $t\in[0, T ]$ with $\triangle t:= \max_{j=1,...,J_n}(t_j- t_{j-1})\rightarrow 0,
n\rightarrow\infty.$
\end{lemma}
\begin{remark}
From the previous lemma we can see that, if the integrand $\phi$
is $\mathcal{F}$-adapted, then the Riemann sums are also an approximation to the It\^{o}
integral of $\phi$ with respect to the Brownian motion. Hence in this case the
forward integral and the It\^{o} integral coincide. In this sense we can regard
the forward integral as an extension of the It\^{o} integral to a nonanticipating
setting.
\end{remark}
In the sequel we give some useful properties of the forward integral. The
following result is an immediate consequence of the definition.

\begin{lemma}
 Suppose $\phi$ is a forward integrable stochastic process and $G$ a
random variable. Then the product $G\phi$ is forward integrable stochastic process and
\begin{equation}
\int_0^TG\phi(t)d^-B(t) = G\int_0^T\phi(t)d^-B(t)
\end{equation}
\end{lemma}
The next result shows that the forward integral is an extension of the
integral with respect to a semimartingale.
\begin{lemma}
Let $\mathbb{G} := \{\mathcal{G}_t, t\in[0, T ]\} (T > 0)$ be a given filtration. Suppose
that
\begin{enumerate}
\item $B$ is a semimartingale with respect to the filtration $\mathbb{G}$.
\item $\phi$ is $\mathbb{G}$-predictable and the integral
\begin{equation}
\int_0^T\phi(t)dB(t),
\end{equation}
with respect to $B$, exists.\\
Then $\phi$ is forward integrable and
\begin{equation}
\int_0^T\phi(t)d^-B(t)=\int_0^T\phi(t)dB(t).
\end{equation}
\end{enumerate}
\end{lemma}
As a consequence of the above we get the following useful result:
\begin{lemma} \label{Lemma7.20}
Let $\varphi(t,y)$ be an $\mathbb{F}$-adapted process for each $y \in \mathbb{R}$ such that
$$\int_0^T \phi(t,y) dB(t)$$
exists for each $y \in \mathbb{R}$. Let $Y$ be a random variable. Then $\varphi(t,Y)$ is forward integrable and
\begin{equation}
\int_0^T \varphi(t,Y)d^{-}B(t) = \int_0^T \varphi(t,y)dB(t)_{y=Y}.
\end{equation}
\end{lemma}

We now turn to the It\^{o} formula for forward integrals. In this connection it is
convenient to introduce a notation that is analogous to the classical notation
for It\^{o} processes.
\begin{definition}
A forward process (with respect to $B$) is a stochastic process
of the form
\begin{equation}\label{forward form 1}
X(t) = x +\int_0^tu(s)ds +\int_0^tv(s)d^-B(s), \quad t\in[0, T ],
\end{equation}
($x$ constant), where $\int_0^T|u(s)|ds <\infty, \mathbf{P}$-a.s.
and $v$ is a forward integrable stochastic process. A shorthand notation for
(\ref{forward form 1}) is that
\begin{equation}
d^-X(t) = u(t)dt + v(t)d^-B(t).
\end{equation}
\end{definition}
\begin{theorem}{The one-dimensional It\^{o} formula for forward integrals.} \\
Let
\begin{equation}
d^-X(t) = u(t)dt + v(t)d^-B(t)
\end{equation}
be a forward process. Let $f\in\mathbf{C}^{1,2}([0, T ]\times\mathbb{R})$ and define
\begin{equation}
Y (t) = f(t,X(t)), \quad t\in[0, T ].
\end{equation}
$Y (t), t\in[0, T ]$, is a forward process and
\begin{equation}
d^-Y (t) = \frac{\partial f}{\partial t}(t,X(t))dt+\frac{\partial f}{\partial x}(t,X(t))d^-X(t)+\frac{1}{2}\frac{\partial^2f}{\partial x^2} (t,X(t))v^2(t)dt.
\end{equation}
\end{theorem}

\section{White Noise and Hida-Malliavin Calculus for $\tilde{N} (.)$}
The construction of a white noise theory and a stochastic derivative (Hida-Malliavin derivative) in the pure jump martingale case follows
the same lines as in the Brownian motion case. In this case, the corresponding Wiener-It\^{o} chaos expansion
Theorem states that any $F \in \mathbf{L}^2(\mathcal{F}_T , \mathbf{P})$ (where, in this case, $\mathcal{F}_t = \mathcal{F}^{(\tilde{N}
)}_t$ is the $\sigma$- algebra generated by $\eta(s)=\int_0^s\int_{\mathbb{R}_0} \zeta \tilde{N} (dr, d\zeta); 0 \leq s \leq t$) can be written as
\begin{equation}\label{}
    F=\sum_{n=0}^{\infty}I_n(f_n); \quad f_n\in\hat{\mathbf{L}}^2((\lambda\times\nu)^n),
\end{equation}
where $\hat{\mathbf{L}}^2((\lambda\times\nu)^n)$ is the space of functions $f_n(t_1,\zeta_1,...,t_n,\zeta_n), t_i\in[0,T],  \zeta_i\in\mathbb{R}_0$ such that $f_n\in \mathbf{L}^2((\lambda\times\nu)^n)$ and $f_n$ is symmetric with respect to the pairs of variables $(t_1, \zeta_1), ..., (t_n, \zeta_n)$.
It is important to note that in this case the $n$-times iterated integral $I_n(f_n)$ is taken with respect to
$\tilde{N}(dt, d\zeta)$ and not with respect to $d\eta(t)$. Thus, we define
\begin{equation}\label{}
    I_n(f_n)=n!\int_0^T\int_{\mathbb{R}_0}\int_0^{t_n}\int_{\mathbb{R}_0}...\int_0^{t_2}\int_{\mathbb{R}_0} f_n(t_1,\zeta_1,...,t_n,\zeta_n)\tilde{N}(dt_1, d\zeta_1)
\end{equation}
for $f_n \in \mathbf{L}^2((\lambda\times\nu)^n)$
The It\^{o} isometry for stochastic integrals with respect to $\tilde{N} (dt, d\zeta)$ then gives the following isometry for
the chaos expansion:
\begin{equation}\label{}
    \|F\|_{\mathbf{L}^2(\mathbf{P})}=\sum_{n=0}^{\infty}n!\|f_n\|^2_{\mathbf{L}^2((\lambda\times\nu)^n)}
\end{equation}
As in the Brownian motion case, we use the chaos expansion to define the Hida-Malliavin derivative. Note that
in this case, there are two parameters $t,\zeta$, where $t$ represents time and $\zeta\neq 0$ represents a generic jump size.

 \subsection{The white noise probability space}
From now on we assume that for every $\epsilon >0$ there exists $\rho > 0$ such that
\begin{equation}\label{eq8.4}
\int_{\mathbb{R} \setminus (-\epsilon,\epsilon)} e^{\rho  \zeta} d\nu(\zeta) < \infty.
\end{equation}

This condition implies that the polynomials are dense in $L^2(\mu)$, where $d\mu(\zeta)=\zeta^2 d\nu(\zeta)$. It also guarantees that the measure $\nu$ integrates all polynomials of degree $\geq 2$.\\

 As in the Brownian motion case, the sample space considered is
$\Omega = \mathbf{S}'(\mathbb{R})$, the space of tempered distributions on $\mathbb{R}$, which is a topological
space. We equip this space with the corresponding Borel $\sigma$-algebra
$\mathcal{F} = \mathcal{B}(\mathbb{R})$. By the Bochner-Minlos-Sazonov theorem,
there exists a probability measure $\mathbf{P}$ such that
\begin{equation}\label{}
    \int_{\Omega}e^{i\langle \omega , f\rangle}\mathbf{P}(d\omega)=\exp(\int_{\mathbb{R}}\psi(f(x))dx); \quad f\in\mathcal{S}(\mathbb{R})
\end{equation}
where
\begin{equation}\label{}
    \psi(\omega)=\int_{\mathbb{R}}(e^{i\omega z}-1-i\omega z)\nu(dz); \quad i=\sqrt{-1}
\end{equation}
and $\langle \omega , f\rangle$ denotes the action of $\omega\in\mathcal{S}'(\mathbb{R})$ on  $f\in\mathcal{S}(\mathbb{R})$.
The triple $(\Omega,\mathcal{F}, \mathbb{P})$ defined above is called the (pure jump) L\'{e}vy white
noise probability space.

Let ${p_j(z)}_{j\geq1}$ defined as in \cite{DOP} then it's an orthogonal basis for $\mathbf{L}^2(\nu)$.
Let ${e_i(t)}_{i\geq1}$ be the Hermite functions. Define
\begin{equation}\label{}
    \delta_{\kappa(i,j)}(t,z)=e_i(t)p_j(z),
\end{equation}
where
\begin{equation}\label{kappa}
    \kappa(i,j)=j+(i+j-2)(i+j-1)/2.
\end{equation}

If $\alpha\in\mathcal{J}$ with Index$(\alpha) = j$ and $|\alpha| = m$, we define $\delta^{\otimes \alpha}$ by\\
$\delta^{\otimes \alpha}(t_1,z_1,...,t_m,z_m)\\
=\delta_1^{\otimes \alpha_1}\otimes...\otimes\delta_j^{\otimes \alpha_j}(t_1,z_1,...,t_m,z_m)\\
=\delta_1(t_1,z_1)...\delta_1(t_{\alpha_1},z_{\alpha_1})...\delta_j(t_{m-\alpha_j+1},z_{m-\alpha_j+1})...\delta_j(t_m,z_m).\\$
We set $\delta_i^{\otimes 0} = 1.$ Finally, we let $\delta^{\hat{\otimes} \alpha}$  denote the symmetrized tensor product
of the $\delta^{\otimes \alpha}$:
For $\alpha\in \mathcal{J}$ define
\begin{equation}\label{}
    K_{\alpha}:=I_{|\alpha|}(\delta^{\hat{\otimes} \alpha})
\end{equation}

As in the Brownian motion case, one can prove that $\{K_{\alpha}\}_{\alpha\in\mathcal{J}}$ are orthogonal
in $\mathbf{L}^2(\mathbf{P})$ and
\begin{equation}\label{}
    \|K_{\alpha}\|^2_{\mathbf{L}^2(\mathbf{P})}=\alpha!.
\end{equation}
We have that, if $|\alpha| = m$,
\begin{equation}\label{}
    \alpha!=\|K_{\alpha}\|^2_{\mathbf{L}^2(\mathbf{P})}=m!\|\delta^{\hat{\otimes} \alpha}\|^2_{\mathbf{L}^2((\lambda\times\nu)^n)}
\end{equation}
 By our construction of $\delta^{\hat{\otimes} \alpha}$  we know that any $f\in\mathcal{L}^2((\lambda,\nu)^n)$ can be written as
\begin{equation}\label{eq6.40}
    f(t_1,z_1,...,t_m,z_m)=\delta_1^{\hat{\otimes} \alpha_1}\otimes...\otimes\delta_j^{\hat{\otimes} \alpha_j}(t_1,z_1,...,t_m,z_m)
\end{equation}
Hence
\begin{equation}\label{}
    I_n(f_n)=\sum_{|\alpha|=n}c_{\alpha}K_{\alpha}
\end{equation}

This gives the following theorem.

\begin{theorem} Chaos expansion.
Any $ F \in \mathbf{L}^2(\mathbf{P})$ has a unique expansion of the form
\begin{equation}
    F=\sum_{\alpha\in\mathcal{J}}c_{\alpha}K_{\alpha}
\end{equation}
with $c_{\alpha} \in \mathbb{R} $. Moreover,
\begin{equation}
    \|F\|^2_{\mathbf{L}^2(\mathbf{P})}=\sum_{\alpha\in\mathcal{J}}\alpha!c^2_{\alpha}
\end{equation}
\end{theorem}
\begin{definition} {The Hida-Kondratiev spaces, $0 \leq \rho \leq 1$.}
\begin{enumerate}
\item {Stochastic test functions $(\mathcal{S})_{\rho}$}. \\
Let $(\mathcal{S})_{\rho}$ consist of all $\phi =\sum_{\alpha\in\mathcal{J}}a_{\alpha}K_{\alpha} \in \mathbf{L}^2(\mathbf{P})$ such that
    \begin{equation}
        \|\phi\|^2_{k,\rho}=\sum_{\alpha\in \mathcal{J}}a^2_{\alpha}(\alpha!)^{1+\rho}(2\mathbb{N})^{k\alpha}<\infty, \text{ for all } k\in\mathbb{N},
    \end{equation}
    equipped with the projective topology, where as before
    \begin{equation}
     (2\mathbb{N})^{k\alpha}=\prod _{j\geq1} (2j)^{k\alpha_j}
    \end{equation}
    if $\alpha=(\alpha_1, \alpha_2,...)\in \mathcal{J}.$

\item {Stochastic distributions $(\mathcal{S})_{-\rho}$}. \\
Let $(\mathcal{S})_{-\rho}$ consist of all expansions
$F =\sum_{\alpha\in\mathcal{J}}b_{\alpha}K_{\alpha}$ such that
\begin{equation}\label{}
    \|F\|^2_{-q,\rho}=\sum_{\alpha\in \mathcal{J}}b^2_{\alpha}(\alpha!)^{1-\rho}(2\mathbb{N})^{-q\alpha}<\infty, \text{ for some }  q\in\mathbb{N},
\end{equation}
endowed with the inductive topology.
\item
As in Section 7 we put $(\mathcal{S})_0=(\mathcal{S})$ and $(\mathcal{S})_{-0}=(\mathcal{S})^{\ast}$, called the Hida test function space and the Hida distribution space, respectively.
\item
The space  $(\mathcal{S})_{-\rho}$ is the dual of $(\mathcal{S})_{\rho}$. In particular, $(\mathcal{S})^{\ast}$ is the dual of $(\mathcal{S})$.\\
If $F =\sum_{\alpha\in\mathcal{J}}b_{\alpha}K_{\alpha} \in (\mathcal{S})_{-\rho}$
 and $\phi =\sum_{\alpha\in\mathcal{J}}a_{\alpha}K_{\alpha} \in (\mathcal{S})_{\rho}$
,then the action of
$F$ on $\phi$ is
\begin{equation}\label{}
    \langle F,\phi\rangle=\sum_{\alpha\in\mathcal{J}}a_{\alpha}b_{\alpha}\alpha!.
\end{equation}
\end{enumerate}

\end{definition}
\subsection{The spaces  $\mathcal{G}$ and $\mathcal{G}^{\ast}$ of Smooth and Generalized Random Variables
}
\begin{definition}\label{defi}
\begin{enumerate}
\item
Let $k\in\mathbb{N}_0$. We say that $f =\sum_{m\geq0}I_m f_m \in \mathbf{L}^2(\mathbf{P})$ belongs to the space $\mathcal{G}_k$ if
\begin{equation}\label{}
    \|f\|_{\mathcal{G}_k}^2:=\sum_{m\geq 0}m!\|f_m\|^2_{\mathbf{L}^2((\lambda\otimes\nu)^m)}e^{2km}<\infty.
\end{equation}
We define the space of smooth random variables $\mathcal{G}$ as
\begin{equation}\label{}
    \mathcal{G}=\bigcap_{k\in \mathbb{N}_0}\mathcal{G}_k.
\end{equation}
The space $\mathcal{G}$ is endowed with the projective topology.

\item We say that a formal expansion
\begin{equation}\label{}
    G=\sum_{m\geq0}I_m(g_m)
\end{equation}
belongs to the space $\mathcal{G}_{-q} (q\in\mathbb{N}_0)$ if
\begin{equation}\label{}
  \|G\|_{\mathcal{G}^2_{-q}}=  \sum_{m\geq0}m!\|g_m\|^2_{\mathbf{L}^2((\lambda\times\nu)^m)}
\end{equation}
The space of generalized random variables $\mathcal{G}^{\ast}$ is defined as
\begin{equation}\label{}
    \mathcal{G}^{\ast}=\bigcup_{q\in\mathbb{N}_0}\mathcal{G}_{-q}
\end{equation}

We equip $\mathcal{G}^{\ast}$ with the inductive topology. Note that $\mathcal{G}^{\ast}$ is the dual of $\mathcal{G}$,
with action
\begin{equation}\label{}
    \langle G, f\rangle= \sum_{m\geq0}m!(f_m,g_m)_{\mathbf{L}^2((\lambda\times\nu)^m)}
\end{equation}

if $G\in \mathcal{G}^\ast$ and $f \in \mathcal{G}$.
Also note that the connection between the expansions
\begin{equation}\label{}
    F=\sum_{m\geq0}I_m(f_m)
\end{equation}
and
\begin{equation}\label{}
   F=\sum_{\alpha\in\mathcal{J}}c_{\alpha}K_{\alpha}
\end{equation}

is given by
\begin{equation}\label{}
    f_m=\sum_{|\alpha|=m}c_{\alpha}\delta ^{\hat{\otimes}\alpha}
\end{equation}

with the functions $\delta ^{\hat{\otimes}\alpha}$ as in (\ref{eq6.40}). Since this gives

\begin{equation}\label{}
    \|I_m(f_m)\|=m!\|f_m\|^2_{\mathbf{L}^2((\lambda\times\nu)^m)}=\sum_{|\alpha|=m}c^2_{\alpha}\|K_{\alpha}\|^2_{\mathbf{L}^2(\mathbf{P})}
\end{equation}

it follows that we can express the $\mathcal{G}_r$-norm of $F$ as follows:
\begin{equation}\label{}
    \|F\|^2_{\mathcal{G}_r}=\sum_{m\geq0}(\sum_{|\alpha|=m}c^2_{\alpha}\|K_{\alpha}\|^2_{\mathbf{L}^2_{(\mathbf{P})}})e^{2rm}; \quad r \in\mathbb{Z}
\end{equation}

By inspecting Theorem(\ref{defi}), we find the following chain of continuous inclusions
\begin{equation}\label{}
    (\mathcal{S})\subset\mathcal{G}\subset\mathbf{L}^2(\mathbf{P}) \subset\mathcal{G}^\ast\subset (\mathcal{S})^{\ast}
\end{equation}
\end{enumerate}
\end{definition}
\subsection{The Hida-Malliavin Derivative on $\mathcal{G}^{\ast}$}
\begin{definition}
Let $F = \sum_{\alpha\in\mathcal{J}}c_{\alpha}K_{\alpha}\in \mathcal{G}^{\ast}$
 Then define the stochastic derivative of $F$ at $(t, z)$ by
 \begin{eqnarray}
   D_{t,z}F &=& \sum_{\alpha}c_{\alpha}\sum_i\alpha_iK_{\alpha-\epsilon^i}.\delta^{\hat{\otimes}\epsilon^i} \\
    &=& \sum_{\alpha}c_{\alpha}\sum_{k,m}\alpha_{\kappa(k,m)}K_{\alpha-\epsilon^{\kappa(k,m)}}.e_k(t)p_m(z) \\ \nonumber
    &=& \sum_{\beta}(\sum_{k,m}c_{\beta+\epsilon^{\kappa(k,m)}}(\beta_{\kappa(k,m)}+1)e_k(t)p_m(z))K_{\beta} \nonumber,
 \end{eqnarray}
 with the map $k(i, j)$ in (\ref{kappa}) and $\epsilon^{l} = \epsilon^{(l)}=(0,0,...,1,0,..,0)$ with $1$ on the l th position .
\end{definition}
We need the following result.
\begin{theorem}
\begin{enumerate}
\item Let $F \in\mathcal{G}^{\ast}$. Then $D_{t,z}F \in\mathcal{G}^{\ast} ,\lambda\times\nu$ - a.e.
\item Let $F\in \mathbf{L}^2(\mathbf{P})$ be $\mathcal{F}_T$ measurable , then $\mathbb{E}[D_{t,z}F|\mathcal{F}_t], t\in[0,T], z\in\mathbb{R}_0$ is an element in $\mathbf{L}^2(\lambda\times\nu\times\mathbf{P})$ and
    \begin{equation}\label{}
        F=\mathbb{E}[F]+\int_0^T\int_{\mathbb{R}_0}\mathbb{E}[D_{t,z}F|\mathcal{F}_t]\tilde{N}(ds,dz).
    \end{equation}
\item Let $F\in \mathbf{L}^2(\mathbf{P})$ then
\begin{equation}\label{}
    \mathbb{E}[\int_0^T\int_{\mathbb{R}_0}F\tilde{N}(ds,dz)]=\mathbb{E}[\int_0^T\int_{\mathbb{R}_0}\mathbb{E}[D_{t,z}F|\mathcal{F}_t]\nu(dz)dt]
\end{equation}
\end{enumerate}
\end{theorem}
\subsection{Ordinary Chain Rules for the Hida-Malliavin Derivative}
\begin{theorem}
\begin{enumerate}
\item \textbf{Product rule.}
 Let $F,G\in\mathbf{D}^{\epsilon}_{1,2}$. Then $FG\in\mathbf{D}^{\epsilon}_{1,2}$ and
 \begin{equation}
      D_{t,z}(F G)= FD_{t,z}G + GD_{t,z}F + D_{t,z}FD_{t,z}G
 \end{equation}
 \item Let $F \in\mathbf{D}^{\epsilon}_{1,2}$
 \begin{equation}
 D_{t,z}(F^n) = (F + D_{t,z}F)^n - F^n.
 \end{equation}
\item Let $F\in\mathbf{D}_{1,2}$ and let $\phi$ be a real continuous
function on $\mathbb{R}$. Suppose $\phi(F)\in\mathbf{L}^2(\mathbf{P})$ and $\phi(F + D_{t,z}F)\in\mathbf{L}^2(\mathbf{P}\times\lambda\times\nu).$
Then $\phi(F) \in \mathbf{D}_{1,2}$ and
\begin{equation}
     D_{t,z}\phi(F) = \phi(F + D_{t,z}F)-\phi(F).
\end{equation}
\end{enumerate}
\end{theorem}

\subsection{The Wick Product}
Now we use the chaos expansion in terms of $\{K_{\alpha}\}_{\alpha\in\mathcal{J}}$ to define the (Poisson random measure type)  Wick product and study some of
its properties.
\begin{definition}
Let $F=\sum_{\alpha\in\mathcal{J}}a_{\alpha}K_{\alpha}$ and $G =\sum_{\beta\in\mathcal{J}}b_{\beta}K_{\beta}$ be two elements of $(\mathcal{S})_{-1}$. Then we define the Wick product of $F$ and $G$ by
\begin{equation*}
    F\diamond G=\sum_{\alpha,\beta\in\mathcal{J}}a_{\alpha}b_{\beta}K_{\alpha+\beta}=\sum_{\gamma\in\mathcal{J}}(\sum_{\alpha+\beta=\gamma}a_{\alpha}b_{\beta})K_{\gamma}
\end{equation*}
\end{definition}
We list some properties of the Wick product:
\begin{enumerate}
\item $F,G \in(\mathcal{S})_1 \Rightarrow F \diamond G \in(\mathcal{S})_1$.
\item $F,G \in(\mathcal{S})_{-1} \Rightarrow F \diamond G \in(\mathcal{S})_{-1}$.
\item $F,G\in(\mathcal{S})^{\ast} \Rightarrow F \diamond G \in(\mathcal{S})^{\ast}$.
\item $F,G \in(\mathcal{S}) \Rightarrow F \diamond G \in(\mathcal{S})$.
\item $F \diamond G = G \diamond F$.
\item $F \diamond (G \diamond H) = (F \diamond G) \diamond H$.
\item $F \diamond (G + H) = F \diamond G + F \diamond H$.
\item $I_n(f_n) \diamond I_m(g_m) = I_{n+m}(f_n\hat{\otimes}g_m)$
\end{enumerate}
In view of the properties $(1)$ and $(4)$ we can define the Wick powers $X^{\diamond n}$
$(n = 1, 2, ...)$ of $X \in(\mathcal{S})_{-1}$ as
\begin{equation*}
  X^{\diamond n} := X\diamond X \diamond ...\diamond X \text{ (n times) }.
  \end{equation*}

We put $X^{\diamond 0} := 1$. Similarly, we define the Wick exponential $\exp^{\diamond}X$ of $X\in(\mathcal{S})_{-1}$  by
\begin{equation}\label{}
  \exp^{\diamond}X :=\sum_{n=0}^{\infty}\frac{1}{n!}X^{\diamond n}.
  \end{equation}
In view of the aforementioned properties, we
have that
\begin{equation}\label{}
  (X + Y )^{\diamond2} = X^{\diamond2} + 2X \diamond Y + Y^{\diamond2}
  \end{equation}

and also
\begin{equation}\label{}
  \exp^{\diamond}(X + Y ) = \exp^{\diamond}X \diamond \exp^{\diamond}Y,
  \end{equation}
for $X, Y\in\mathcal{S}_{-1}$
As before let $\mathbb{E}[.]$ denote the generalized expectation. Then
we see that
\begin{equation}\label{}
    \mathbb{E}[X \diamond Y ] = \mathbb{E}[X]\mathbb{E}[Y ],
\end{equation}
for $X, Y \in(\mathcal{S})_{-1}$. Note that independence is not required for this identity to
hold.
By induction, it follows that
\begin{equation}\label{}
    \mathbb{E}[\exp^{\diamond}X] = \exp{\mathbb{E}[X]},
\end{equation}
for $X \in(\mathcal{S})_{-1}$.

\begin{example}
Choose $h\in \mathbf{L}^2([0, T ])$ and define $F=\int_0^Th(t)d\eta(t)$. Then
\begin{eqnarray}
  F\diamond F &=& I_1(h(t)z)\diamond I_1(h(t)z)\nonumber \\
   &=& I_2(h(t_1)h(t_2)z_1z_2)\nonumber \\
   &=&2\int_0^T\int_{\mathbb{R}}(\int_0^T \int_{\mathbb{R}}h(t_1)h(t_2)z_1z_2\tilde{N}(dt_1,dz_1))\tilde{N}(dt_2,dz_2)\nonumber \\
   &=& 2 \int_0^T(\int_0^{t_2}h(t_1)d\eta(t_1)h(t_2))h(t_2)d\eta(t_2)
\end{eqnarray}
By the It\^{o} formula, if we put $X(t) :=\int_0^T h(s)dB(s)$,
\begin{equation*}
    d(X(t))^2=2X(t)dX(t)+h^2(t)\int_{\mathbb{R}}z^2N(dt,dz).
\end{equation*}
Hence
\begin{equation}\label{}
    F\diamond F=2\int_0^TX(s)dX(s)=X^2(T)-\int_0^T\int_{\mathbb{R}}h^2(s)z^2N(ds,dz).
\end{equation}
In particular, choosing $h = 1$ we get
\begin{equation}\label{}
    \eta(T)\diamond\eta(T)=\eta^2(T)-\int_0^T\int_{\mathbb{R}}z^2N(ds,dz).
\end{equation}
\end{example}
\begin{example} {Wick/Dol\'{e}ans-Dade exponential.} \\
Choose $\gamma\geq-1$ deterministic
such that
\begin{equation}\label{}
    \int_{\mathbb{R}}\int_{\mathbb{R}}\{|\ln(1+\gamma(t,z))|+\gamma^2(t,z)\}\nu(dz)dt<\infty
\end{equation}
and put
\begin{equation}\label{}
   F= \exp^{\diamond}(\int_{\mathbb{R}}\int_{\mathbb{R}_0}\gamma(t,z)\tilde{N}(dt,dz))
\end{equation}
To find an expression for $F$ not involving the Wick product, we proceed as
follows:
Define
\begin{eqnarray}
  Y(t) &=& \exp^{\diamond}(\int_0^t\int_{\mathbb{R}_0}\gamma(s,z)\tilde{N}(ds,dz)) \nonumber\\
   &=& \exp^{\diamond}(\int_0^t\int_{\mathbb{R}_0}\gamma(s,z) \dot{\tilde{N}}(s,z)\nu(dz)ds)
\end{eqnarray}
Then we have
\begin{equation}\label{}
    \frac{dY(t)}{dt}=Y(t)\diamond\int_{\mathbb{R}_0}\gamma(t,z) \dot{\tilde{N}}(t,z)\nu(dz)
\end{equation}
or
\begin{equation}\label{}
    dY(t)=Y(t)\int_{\mathbb{R}_0}\gamma(t,z)\tilde{N}(dt,dz)
\end{equation}
Using It\^{o} calculus the solution of this equation is
\begin{eqnarray}
  Y(t) &=&Y(0)\exp(\int_0^t\int_{\mathbb{R}_0}\{\ln(1+\gamma(s,z))-\gamma(s,z)\}\nu(dz)ds  \\
   &+& \int_0^t\int_{\mathbb{R}_0}\ln(1+\gamma(s,z)) \tilde{N}(ds,dz))
\end{eqnarray}
Comparing the two expressions for $Y (t)$, we conclude that
\begin{eqnarray}
   \exp^{\diamond}(\int_{\mathbb{R}}\int_{\mathbb{R}_0}\gamma(s,z)\tilde{N}(ds,dz))&=&\exp(\int_{\mathbb{R}}\int_{\mathbb{R}_0}\{\ln(1+\gamma(s,z))-\gamma(s,z)\}\nu(dz)ds  \nonumber\\
   &+& \int_0^t\int_{\mathbb{R}_0}\ln(1+\gamma(s,z)) \tilde{N}(ds,dz))
\end{eqnarray}

In particular, choosing
\begin{equation}\label{}
    \gamma(s,z)=h(s)z \mbox{  with } h\in\mathbf{L}^2(\mathbb{R})
\end{equation}
we get
\begin{eqnarray}
  \exp^{\diamond}(\int_{\mathbb{R}}h(s)d\eta(s)) &=& \exp(\int_{\mathbb{R}}\int_{\mathbb{R}_0}\{\ln(1+h(s)z)-h(s)z\}\nu(dz)ds \nonumber \\
  &+&  \int_{\mathbb{R}}\int_{\mathbb{R}_0}\ln(1+h(s)z)\tilde{N}(ds,dz))
\end{eqnarray}
\end{example}
\subsection{The Wick chain rule for a Poisson random measure}
\begin{theorem}
Let $F\in (\mathcal{S})_{-1}$ and $\phi\in\mathbf{C}^1(\mathbb{R})$ then:
\begin{equation}
  D_{t,z} \phi^{\diamond}(F) = (\phi ' )^{\diamond}(F) \diamond D_{t,z}F
\end{equation}
\end{theorem}
\dproof
First note that if $ \psi(s,\zeta)$ is deterministic then
\begin{align}
\exp ^{\diamond} \big( \int_0^T\int_{\mathbb{R}} \psi(s,\zeta) \tilde{N}(ds,d\zeta)\big)
=K \exp \big( \int_0^T\int_{\mathbb{R}} \ln (1+\psi(s,\zeta))\tilde{N}(ds,d\zeta) \big)
\end{align}
where
\begin{equation}
K:= \exp \big( \int_0^T \int_{\mathbb{R}} [\ln (1+\psi(s,\zeta)) -\psi(s,\zeta)]\nu(d\zeta)ds \big)
\end{equation}
Using this and the Wick rule for $\exp $ we get
\begin{align}
&D_{t,z} \exp^{\diamond} \big( \int_0^T\int_{\mathbb{R}} \psi(s,\zeta)\tilde{N}(ds,d\zeta) \big)\nonumber\\
& = K D_{t,z} \exp \big( \int_0^T\int_{\mathbb{R}} \ln(1+\psi(s,\zeta))\tilde{N}(ds,d\zeta) \big)  \nonumber\\
&=K \big[ \exp \big( \int_0^T\int_{\mathbb{R}} \ln(1+\psi(s,\zeta))\tilde{N}(ds,d\zeta) + \ln(1+\psi(t,z))\big) \nonumber\\
&-\exp(\int_0^T\int_{\mathbb{R}} \ln(1+\psi(s,\zeta))\tilde{N}(ds,d\zeta))\big ] \nonumber\\
&=K  \exp \big( \int_0^T\int_{\mathbb{R}} \ln(1+\psi(s,\zeta))\tilde{N}(ds,d\zeta)\big[\exp( \ln(1+\psi(t,z)))-1\big]\big) \nonumber\\
&= \exp^{\diamond} \big( \int_0^T\int_{\mathbb{R}} \psi(s,\zeta)\tilde{N}(ds,d\zeta) \big) \psi(t,z)
\end{align}
In view of this, and the fact that the set of linear combinations of exponentials of the form above are dense in $L^2(P)$ we have the result.
\fproof

\subsection {Conditional expectation}

If $F= \sum_{\alpha\in\mathcal{J}}b_{\alpha}K_{\alpha} \in (\mathcal{S})_{-1},$ we define its conditional expectation by the expansion
 \begin{equation}
\mathbb{E}[F|\mathcal{F}_t] :=  \sum_{\alpha\in\mathcal{J}}b_{\alpha}\mathbb{E}[K_{\alpha}| \mathcal{F}_t]
\end{equation}
Then as in Section 7 we have::
\begin{enumerate}
\item
If $F,G \in (\mathcal{S})_{-1},$ then $\mathbb{E}[(F \diamond G) | \mathcal{F}_t] \in (\mathcal{S})_{-1}$ and \\
$\mathbb{E}[(F \diamond G) | \mathcal{F}_t] = \mathbb{E}[F|\mathcal{F}_t] \diamond \mathbb{E}[G|\mathcal{F}_t].$
\item
If $F,G \in (\mathcal{G})^{\ast},$ then $\mathbb{E}[(F \diamond G) | \mathcal{F}_t] \in (\mathcal{G})^{\ast}.$
\end{enumerate}
\subsection{The Forward Integral for L\'{e}vy processes}
\begin{definition}
The forward integral
\begin{equation}
J(\theta) :=\int_0^T\int_{\mathbb{R}_0}\theta(t, z)\tilde{N}(d^-t, dz)
\end{equation}
with respect to the Poisson random measure $\tilde{N}$ of a stochastic function (or
random field) $\theta(t, z), t\in[0, T ], z\in\mathbb{R}_0 (T >0)$, with
\begin{equation}
\theta(t, z) := \theta(\omega, t, z), \omega\in\Omega,
\end{equation}
and c\`{a}gl\`{a}d with respect to $t$, is defined as
\begin{equation}
J(\theta) = \lim_{m\rightarrow\infty}\int_0^T\int_{\mathbb{R}_0}\theta(t, z)\mathbf{1}_{U_m}\tilde{N}(dt, dz)
\end{equation}
if the limit exists in $\mathbf{L}^2(\mathbf{P})$. Here, $U_m, m = 1, 2, ...,$ is an increasing sequence of
compact sets $U_m \subseteq\mathbb{R}_0:= \mathbb{R}\setminus \{0\}$ with $\nu(U_m) < \infty$ such that $\lim_{m\rightarrow\infty}U_m=\mathbb{R}_0$.
\end{definition}
\begin{remark}
Note that if $\mathbb{H} := \{\mathcal{H}_t, t\in[0, T ]\}$ is a filtration such that
\begin{enumerate}
\item $\mathcal{F}_t\subseteq \mathcal{H}_t$ for all $t$
\item The process $\eta(t) =\int_0^t\int_{\mathbb{R}_0}z \tilde{N}(ds, dz), t\in[0, T ]$, is a semimartingale with
respect to $\mathbb{H}$
\item The stochastic process $\theta= \theta(t, z), t\in[0, T ], \in\mathbb{R}_0$ is $\mathbb{H}$-predictable and
\item The integral $\int_0^t\int_{\mathbb{R}_0}\theta(t, z)\tilde{N}(dt, dz)$ exists as a classical It\^{o} integral
\end{enumerate}
then the forward integral of $\theta$ with respect to $\tilde{N}$ also exists and we have
\begin{equation}
\int_0^T\int_{\mathbb{R}_0}\theta(t, z)\tilde{N}(d^-t, dz)=\int_0^T\int_{\mathbb{R}_0}\theta(t, z)\tilde{N}(dt, dz)
\end{equation}
This follows from the basic construction of the semimartingale integral. Thus, the forward integral can be regarded as an extension of the
It\^{o} integral to possibly nonsemimartingale contexts.
\end{remark}
\begin{remark}
Directly from the definition we can see that if $G$ is a random
variable then
\begin{equation}
G\int_0^T\int_{\mathbb{R}_0}\theta(t, z)\tilde{N}(d^-t, dz)=\int_0^T\int_{\mathbb{R}_0}G\theta(t, z)\tilde{N}(d^-t, dz)
\end{equation}
\end{remark}
As a consequence of the above we get the following useful result (compare with Lemma 7.20):
\begin{lemma} \label{Lemma8.12}
Let $\varphi(t,y,\zeta)$ be an $\mathbb{F}$-adapted process for each $y \in \mathbb{R}$ such that
$$\int_0^T \phi(t,y,\zeta) \tilde{N}(dt,d\zeta)$$
exists for each $y \in \mathbb{R}$. Let $Y$ be a random variable. Then $\varphi(t,Y,\zeta)$ is forward integrable and
\begin{equation}
\int_0^T \varphi(t,Y,\zeta)\tilde{N}(d^{-}t,\zeta) = \int_0^T \varphi(t,y,\zeta)\tilde{N}(dt,d\zeta)_{y=Y}.
\end{equation}
\end{lemma}

\begin{definition}
A forward process is a measurable stochastic function
$X(t) = X(t,\omega), t\in[0, T ], \omega\in\Omega, $ that admits the representation
\begin{equation}\label{forward form}
X(t) = x +\int^t_0\int_{\mathbb{R}_0}\theta(s, z)\tilde{N}(d^-s, dz) +\int^t_0\alpha(s)ds,
\end{equation}
where x = X(0) is a constant. A shorthand notation for (\ref{forward form}) is
\begin{equation}
d^-X(t) =\int_{\mathbb{R}_0}\theta(t, z) \tilde{N}(d^-t, dz) + \alpha(t)dt, X(0) = x.
\end{equation}
We call $d^-X(t)$ the forward differential of $X(t), t\in[0, T ]$.
\end{definition}
\begin{theorem}It\^{o} formula for forward integrals.\\
 Let $X(t), t\in[0, T ]$,
be a forward process of the form $(\ref{forward form})$, where $\theta(t, z), t\in[0, T ], z\in\mathbb{R}_0$, is
locally bounded in $z$ near $z = 0,  \mathbf{P}\times\lambda$- a.e., such that
\begin{equation}
\int_0^T\int_{\mathbb{R}_0}|\theta(t, z)|^2\nu(dz)dt <\infty \quad a.s. \mathbf{P}.
\end{equation}
Also suppose that $|\theta(t, z)|, t\in[0, T ], z\in\mathbb{R}_0$, is forward integrable. For any
function $f\in\mathbf{C}^2(\mathbb{R})$, the forward differential of $Y (t) = f(X(t)), t\in[0, T ]$, is
given by the following formula:
\begin{align}
&d^-Y (t) = f'(X(t))\alpha(t)dt+\int_{\mathbb{R}_0}\big(f(X(t^-) + \theta(t, z))-f(X(t^-))-f'(X(t^-))\theta(t, z)\big)\nu(dz)dt\nonumber\\
&+\int_{\mathbb{R}_0}\big(f(X(t^-) + \theta(t, z))-f(X(t^-))\big)\tilde{N}(d^-t, dz).
\end{align}
\end{theorem}


\begin{thebibliography}{9999}

\bibitem[A\O ]{AO}
N. Agram and B. \O ksendal:
Malliavin calculus and optimal control of stochastic Volterra equations.
J. Optim. Theory Appl. (2015) DO1 10. 1007/s 10957-015-0753-5.

\bibitem[Aa\O PU]{AaOPU}
K. Aase, B. \O{}ksendal, N. Privault and J. Ub\o e: White noise generalizations of the Clark-Haussmann-Ocone theorem with application to mathematical finance. Finance Stoch. 4 (2000), 465-496.

\bibitem[Aa\O U]{AaOU}
K. Aase, B. \O{}ksendal and J. Ub\o e: Using the Donsker delta function to compute hedging strategies. Potential Analysis 14 (2001), 351-374.

\bibitem[B]{B} L. Breiman: Probability. Addison-Wesley 1968.

\bibitem[BBS]{BBS}
O. E. Barndorff-Nielsen, F.E. Benth and B. Szozda: On stochastic integration for volatility modulated Brownian-driven Volterra processes via white noise analysis. arXiv:1303.4625v1, 19 March 2013.

\bibitem[B\O]{BO}
F. Biagini and B. \O{}ksendal: A general stochastic calculus approach to insider trading.Appl. Math. \& Optim. 52 (2005), 167-181.

\bibitem [DM\O R]{DMOR}K.R. Dahl, S.-E. A. Mohammed, B. \O ksendal and E. R. R\o se; Optimal control with noisy memory and BSDEs with Malliavin derivatives. arXiv: 1403.4034 (2014).

\bibitem[DM\O P1]{DMOP1} G. Di Nunno, T. Meyer-Brandis, B. {\O}ksendal and F. Proske: Malliavin calculus and anticipative It\^ o formulae for L\'{e}vy processes. Inf. Dim. Anal. Quantum Prob. Rel. Topics 8 (2005), 235-258.

\bibitem[DM\O P2]{DMOP2} G. Di Nunno, T. Meyer-Brandis, B. {\O}ksendal and F. Proske: Optimal portfolio for an insider in a market driven by L\'{e}vy processes. Quant. Finance 6 (2006), 83-94.

\bibitem[D\O ]{DiO} G. Di Nunno and B. {\O}ksendal: The Donsker delta function, a representation formula for functionals of a L\'{e}vy process and application to hedging in incomplete markets. S\'{e}minaires et Congr\`{e}es, Societ\'{e} Math\'{e}matique de France, Vol. 16 (2007), 71-82.

\bibitem[D\O ]{DO} G. Di Nunno and B. \O ksendal: A representation theorem and a sensitivity result for functionals of jump diffusions. In A.B. Cruzeiro, H. Ouerdiane and N. Obata (editors): Mathematical Analysis and Random Phenomena. World Scientific 2007, pp. 177 - 190.

\bibitem[D\O P]{DOP} G. Di Nunno, B. {\O}ksendal and F. Proske: Malliavin Calculus for L\'{e}vy Processes with Applications to Finance. Universitext, Springer 2009.

\bibitem[H\O UZ]{HOUZ} H. Holden, B. {\O}ksendal, J. Ub\o e and T. Zhang: Stochastic Partial Differential Equations. Universitext, Springer, Second Edition 2010.

\bibitem[LP]{LP} A. Lanconelli and F. Proske: On explicit strong solution of It\^ o-SDEs and the Donsker delta function of a diffusion.
Inf. Dim. Anal. Quatum Prob Rel. Topics 7 (2004),437-447.

\bibitem[M\O P]{MOP} S. Mataramvura, B. \O ksendal and F. Proske: The Donsker delta function of a L\'{e}vy process with application to chaos expansion of local time. Ann. Inst H. Poincar\'{e} Prob. Statist. 40 (2004), 553-567.

\bibitem[MP]{MP} T. Meyer-Brandis and F. Proske: On the existence and explicit representability of strong solutions of L\'{e}vy noise driven SDEs with irregular coefficients. Commun. Math. Sci. 4 (2006), 129-154.

\bibitem[\O R1]{OR1} B. \O ksendal and E. R\o se: A white noise approach to insider trading. Manuscript University of Oslo, 20 November 2014

\bibitem[\O R2]{OR2} B. \O ksendal and E. R\o se: Applications of white noise to mathematical finance. Manuscript University of Oslo, 5 February 2015

\bibitem[\O S1]{OS1}
B. \O{}ksendal and A. Sulem: Applied Stochastic Control of Jump Diffusions. Second Edition. Springer 2007

\bibitem[\O S2]{OS2}
B. \O{}ksendal and A. Sulem: Risk minimization in financial markets modeled by It\^ o-L\' evy processes. Afrika Matematika (2014), DOI: 10.1007/s13370-014-02489-9.

\bibitem[P]{P} P. Protter: Stochastic Integration and Differential Equations. Second Edition. Springer 2005

\bibitem[PK]{PK} I. Pikovsky and I. Karatzas: Anticipative portfolio optimization. Adv. Appl. Probab. 28 (1996), 1095-1122.

\bibitem[RV]{RV} F. Russo and P. Vallois: Forward, backward and symmetric stochastic integration. Probab. Theor. Rel. Fields 93 (1993), 403-421.

\bibitem[RV1]{RV1} F. Russo and P. Vallois. The generalized covariation process and It\^{o} formula. Stoch. Proc. Appl., 59(4):81-104, 1995.

\bibitem[RV2]{RV2} F. Russo and P. Vallois. Stochastic calculus with respect to continuous finite quadratic variation processes. Stoch. Stoch. Rep., 70(4):1-40, 2000.


\end{thebibliography}
\end{document}